\RequirePackage{fix-cm}
\documentclass[smallextended]{svjour3}       
\smartqed  
\usepackage{graphicx}

\RequirePackage{amsmath}
\RequirePackage[numbers]{natbib}
\RequirePackage[colorlinks,citecolor=blue,urlcolor=blue]{hyperref}

\usepackage{amssymb}
\usepackage{algorithm}
\usepackage{algorithmic}
\usepackage{color}
\usepackage{mathtools}
\usepackage{enumitem} 
\usepackage{dsfont} 
\usepackage[toc,page]{appendix}

\def\ba#1\ea{\begin{align*}#1\end{align*}}
\def\ban#1\ean{\begin{align}#1\end{align}}

\newcommand\numberthis{\addtocounter{equation}{1}\tag{\theequation}} 

\usepackage{bbold}

\newcommand{\wi}{\ensuremath{w_i} }
\newcommand{\wtj}{\ensuremath{\widetilde{ w}_{j}} }

\allowdisplaybreaks 
\def\argmax{\operatornamewithlimits{arg\,max}} 

\newcommand{\bE}{\ensuremath{\mathds{E}} }
\newcommand{\bP}{\ensuremath{\mathds{P}} }
\newcommand{\bR}{\ensuremath{\mathds{R}} }

\newcommand{\bEq}[1]{ \mathds{E} \left[\left. #1\right|{\bf q}(t)={\bf q}\right]  }

\newcommand{\vq}{\ensuremath{{\bf q}} }
\newcommand{\bvq}{\ensuremath{\overline{{\bf q}}} } 
\newcommand{\va}{\ensuremath{{\bf a}} }
\newcommand{\bva}{\ensuremath{\overline{{\bf a}}} } 
\newcommand{\vs}{\ensuremath{{\bf s}} }
\newcommand{\vu}{\ensuremath{{\bf u}} }
\newcommand{\valpha}{\ensuremath{{\mbox{\boldmath{$\alpha$}}}} } 
\newcommand{\vlam}{\ensuremath{{\mbox{\boldmath{$\lambda$}}}} } 
\newcommand{\vsig}{\ensuremath{{\mbox{\boldmath{$\sigma$}}}} } 
\newcommand{\vnu}{\ensuremath{{\mbox{\boldmath{$\nu$}}}} } 
\newcommand{\vzeta}{\ensuremath{{\mbox{\boldmath{$\zeta$}}}} } 
\newcommand{\vchi}{\ensuremath{{\mbox{\boldmath{$\chi$}}}} } 
\newcommand{\vk}{\ensuremath{{\mbox{\boldmath{$k$}}}} } 
\newcommand{\vqt}{\ensuremath{{\bf q}(t)} }

\newcommand{\vx}{\ensuremath{{\bf x}} }
\newcommand{\vy}{\ensuremath{{\bf y}} }

\newcommand{\vf}{\ensuremath{{\bf f}} }

\newcommand{\vei}{\ensuremath{{\bf e}^{(i)}} }
\newcommand{\vetj}{\ensuremath{\widetilde{{\bf e}}^{(j)}} }
\newcommand{\ei}{\ensuremath{e^{(i)}} }
\newcommand{\etj}{\ensuremath{\widetilde{ e}^{(j)}} }
\newcommand{\ki}{\ensuremath{\kappa_i }}
\newcommand{\ktj}{\ensuremath{\widetilde{\kappa}_j} }

\newcommand{\vone}{\ensuremath{{\bf 1}} } 

\newcommand{\qij}{\ensuremath{q_{ij}} }
\newcommand{\bqij}{\ensuremath{\overline{q}_{ij}} } 
\newcommand{\qijt}{\ensuremath{q_{ij}(t)} }
\newcommand{\bq}{\ensuremath{\overline{q}} }
\newcommand{\aij}{\ensuremath{a_{ij}} }
\newcommand{\baij}{\ensuremath{\overline{a}_{ij}} } 
\newcommand{\sij}{\ensuremath{s_{ij}} }
\newcommand{\uij}{\ensuremath{u_{ij}} }
\newcommand{\nuij}{\ensuremath{\nu_{ij}} }
\newcommand{\lij}{\ensuremath{_{ij}} }
\newcommand{\para}{\ensuremath{_{\parallel}} }
\newcommand{\per}{\ensuremath{_{\perp}} }
\newcommand{\se}{\ensuremath{^{(\epsilon)}} }
\newcommand{\numin}{\ensuremath{\nu_{\min}} }
\newcommand{\numinp}{\ensuremath{\nu_{\min}'} }

\newcommand{\cC}{\ensuremath{\mathcal{C}} }
\newcommand{\cS}{\ensuremath{\mathcal{S}} }
\newcommand{\cK}{\ensuremath{\mathcal{K}} }
\newcommand{\cF}{\ensuremath{\mathcal{F}} }

\newcommand{\La}{\left \langle }
\newcommand{\Ra}{\right \rangle}

\newcommand{\V}{\ensuremath{V} }
\newcommand{\Vq}{\ensuremath{V(\vq)} }
\newcommand{\Vi}{\ensuremath{V_1} }

\newcommand{\Wper}{\ensuremath{W_{\perp\cK}} }
\newcommand{\Wperq}{\ensuremath{W_{\perp\cK}(\vq)} }

\newcommand{\ri}{\ensuremath{\gamma_i} }
\newcommand{\rtj}{\ensuremath{\widetilde{\gamma}_j} }
\newcommand{\parainc}{\ensuremath{_{\parallel\cK}} }
\newcommand{\perinc}{\ensuremath{_{\perp\cK}} }
\newcommand{\parains}{\ensuremath{_{\parallel\cS}} }
\newcommand{\perins}{\ensuremath{_{\perp\cS}} }

\newcommand{\Biii}{B_3}
\newcommand{\Biv}{B_4}

\newcommand{\hn}{\hspace{-.25em}}

\newcommand{\linc}{\ensuremath{_{n_1n_2}} }
\newcommand{\for}{\ensuremath{\text{ for }} }

\usepackage{subfiles}

\newcommand{\paras}{\ensuremath{_{\parallel \cS}} }
\newcommand{\norm}[1]{\left\lVert#1\right\rVert}

\newcommand{\flij}{\ensuremath{f_{l\lij}} }
\newcommand{\nuli}{\ensuremath{v_{l_i}} } 
\newcommand{\nutlj}{\ensuremath{\widetilde{v}_{l_j}} }
\newcommand{\su}[2]{ \sum_{#1=1}^{n_#2}  }

\newcommand{\w}{\ensuremath{W }}
\newcommand{\wt}{\ensuremath{\widetilde{W}} }
\newcommand{\wk}{\ensuremath{w_k }}
\newcommand{\wtl}{\ensuremath{\widetilde{w}_l} }
\newcommand{\no}{\ensuremath{n_1} }
\newcommand{\nt}{\ensuremath{n_2} }

 \journalname{Queueing Systems}
\begin{document}

\title{Optimal Heavy-Traffic Queue Length Scaling in an Incompletely Saturated Switch}

\titlerunning{Heavy-Traffic Queue Length Scaling in a Switch}        

\author{Siva Theja Maguluri         \and
        Sai Kiran Burle 	\and
        R. Srikant
}


\institute{ Siva Theja Maguluri \at
          Mathematical Sciences \\
          IBM T. J. Watson Research Center \\
          \email{siva.theja@gmail.com}
           \and
           Sai Kiran Burle \at
           Department of ECE and CSL\\
           University of Illinois at Urbana-Champaign\\
           \email{burle2@illinois.edu}
           \and
           R. Srikant \at
           Department of ECE and CSL\\
           University of Illinois at Urbana-Champaign\\
           \email{rsrikant@illinois.edu}
}

\date{Received: date / Accepted: date}

\maketitle

\begin{abstract}
We consider an input queued switch operating under the MaxWeight scheduling algorithm. This system is interesting to study because it is a model for Internet routers and data center networks. Recently, it was shown that the MaxWeight algorithm has optimal heavy-traffic queue length scaling when all ports are uniformly saturated.
Here we consider the case when an arbitrary number of ports are saturated (which we call the incompletely saturated case), and each port is allowed to saturate at a different rate.
We use a recently developed drift technique to show that the heavy-traffic queue length under the MaxWeight scheduling algorithm has optimal scaling with respect to the switch size even in these cases.
\end{abstract}

\section{Introduction}

The $n\times n$ switch is a model that has been widely used and studied to understand the behavior of Internet routers and data center switches. The importance of this model in the design of Internet routers is well known. In the mid to late 1990s when the Internet was exploding, it served as an important model to study and design scheduling algorithms for the switch fabric of Internet routers. The model is now used to understand the design of data centers used for cloud-computing services. Today's data centers consist of a massive number of servers organized in racks, which are interconnected through a data center network.
An ideal data center network is a huge input queued switch with one port for each server. However, real switches are much smaller and they have to be interconnected appropriately, and routing and scheduling algorithms have to be designed, so that the overall network emulates an $n\times n$ switch. Designing such a network is a challenging and active area of research; see \cite{pfabric,fastpass} for example. Here, we do not explicitly consider a data center network, but only consider the $n\times n$ switch which is the underlying model (see \cite{fastpass} which argues why the model is appropriate even for a data center network) and study the behavior of the well-known MaxWeight algorithm \cite{TasEph_92} for this model.

As mentioned in \cite{shah_switch_open}, the $n\times n$ switch model also serves as a canonical theoretical example of a problem which exhibits the so-called multi-dimensional state-space collapse, which makes it difficult to study using traditional heavy-traffic theory. Recently, it has been shown in \cite{switch_arxiv} that the heavy-traffic behavior of the mean queue length in an $n\times n$ switch operating under the MaxWeight scheduling algorithm can be precisely characterized using a non-trivial extension of a drift technique introduced in \cite{erysri-heavytraffic}. In particular, one of the key contributions of \cite{switch_arxiv} is to extend the drift technique to cover the case of multi-dimensional state-space collapse. The result in \cite{switch_arxiv} also resolved an open question on the scaling behavior of the heavy-traffic queue length in a switch operating under the MaxWeight algorithm. In particular, it showed that the total heavy-traffic-scaled queue length is $O(n)$ or the mean heavy-traffic-scaled delay experienced by a packet is $O(1).$ While there have been other results establishing $O(1)$ delay scaling, the significance of the result in \cite{switch_arxiv} is that the result holds for the original MaxWeight algorithm introduced in \cite{TasEph_92}, with no additional scheduling operations required.

The results in \cite{switch_arxiv} were obtained under the assumption that every input and output port of the switch is saturated (i.e., close to capacity), and the arrival rates to each input port and each output port are close to capacity by the same amount. For the purposes of this paper, we call a switch where only some of ports are saturated an ``incompletely saturated switch." The main purpose of this paper is to show that the MaxWeight algorithm has order-optimal scaling (in the number of ports of the switch) for the case of an incompletely saturated switch and for the case where each port has a different rate of saturation.
The results in \cite{switch_arxiv} were obtained by setting the drift of a test function to zero in steady-state. This function was carefully chosen based on the underlying symmetry when all ports are saturated. Since we do not have such a symmetry in the incompletely saturated switch, a new test function has to be used.
The main reason for this is that the geometry of the state-space collapse is different here than in the all-ports saturated case considered in \cite{switch_arxiv}.
In this paper, we  propose a novel test function to use for the incompletely saturated switch. Another major contribution of the paper is that we express this function in a simple form, so that it can be generalized for use in other problems that exhibit state space collapse into other regions.
The case where each port has a different rate of saturation
is of theoretical importance since it corresponds to the situation where the drift vector is not the identity matrix in the diffusion limit \cite{kang2009diffusion}. In fact, the diffusion limit is no longer symmetric in the components of the limiting stochastic process (i.e., the diffusion limit is not symmetric across the ports), but here we show that the technique in \cite{switch_arxiv} works even in this asymmetric situation to produce an exact formula for a certain linear combination of the queue lengths. This result can be further used to show optimal queue length scaling under some conditions on the saturation rates.

The rest of the paper is organized as follows. The model of an $n\times n$ switch and the MaxWeight algorithm are presented in Section~\ref{sec: system model}. General results on queue lengths are presented in Section~\ref{sec:incomplete}. In Section~\ref{sec:discussion}, we present various extensions and special cases in order to interpret the general result.
Concluding remarks are provided in the last section.

We remark that a very preliminary version of this paper appeared in \cite{sigmetrics_switch}.

\subsection{Related Work}

The MaxWeight algorithm for general stochastic networks, of which the $n\times n$ switch is a special case, was presented in \cite{TasEph_92}, where it was shown that the algorithm is throughput-optimal. The special case of the switch was considered in \cite{mckeown96walrand}, where it was shown that the simpler algorithms such as MaxSize and Maximal matchings are not throughput-optimal. The case of non-stochastic arrivals was considered in \cite{weller1997scheduling}, where in addition to the throughput-optimality of MaxWeight-type algorithms, a lower bound on the throughput loss of simpler algorithms such as Maximal matching was established.

Here we are interested in performance metrics beyond throughput optimality. In particular, we are interested in understanding whether the MaxWeight matching algorithm for switches achieves small queue lengths, at least under a heavy-traffic scaling regime. Using diffusion limits, the heavy-traffic optimality of the MaxWeight algorithm in a switch where only one port is saturated was established in \cite{stolyar2004maxweight}, although the final step of interchanging the order of the heavy traffic scaling limit and letting time go to infinity was not undertaken there. Motivated by this result, \cite{erysri-heavytraffic} studied the switch directly in steady-state, established heavy-traffic optimality, and introduced a new drift method of studying stochastic networks in heavy-traffic. However, it should be emphasized that the results in \cite{stolyar2004maxweight,erysri-heavytraffic} apply only to the case of a single saturated port since they both rely on state of the system collapsing to a single dimension in the heavy-traffic limit. In a recent development, the case where all ports are uniformly saturated (thus, leading to the more difficult case of multi-dimensional state-space collapse) was studied in \cite{switch_arxiv}, where an exact expression for the heavy-traffic scaled queue length under the MaxWeight algorithm is derived. Additionally, this expression shows that the algorithm has heavy-traffic optimal scaling in the size of the switch, resolving an open conjecture stated in \cite{shah_switch_open}. The results in \cite{switch_arxiv} use and significantly extend the drift technique presented in \cite{erysri-heavytraffic}.

State collapse in the case where multiple ports are saturated has been established in \cite{Stolyar_cone_SSC,ShaWis_11,singh2015maxweight} and using the state-space collapse result in \cite{ShaWis_11}, a diffusion limit was established in \cite{kang2009diffusion}. However, properties of the diffusion limit (such as its steady-state distribution or mean queue lengths) were difficult to establish. The result in \cite{switch_arxiv} can thus be interpreted as a derivation of the sum of the first moments of the limiting vector stochastic process, but obtained without going through the usual fluid/diffusion limit scaling arguments. An entirely different technique to study heavy-traffic optimality was presented in \cite{shah2012optimal} where the authors approximate the scheduling decisions made by a switch which can change its schedule infinitely often to simulate a queueing network with product-form steady-state distribution as in \cite{bonald2002insensitivity}. The resulting algorithm is heavy-traffic optimal, but has a very high computational complexity. The optimal scaling of the queue length as a function of the switch size in the non-heavy-traffic limit appears to be still open. Alternatively, one can consider asymptotic regimes other than the heavy-traffic limit. The best known results in this regard are the ones in \cite{neely2007logarithmic,zhong2014scaling}, but these require algorithms that are more involved than the original MaxWeight algorithm.

\section{System Model and Background}\label{sec: system model}
In this section, we present the model of an input queued switch, the MaxWeight scheduling algorithm, and some lemmas that will later play a key role in the results.

\emph{Note on Notation:}
For ease of understanding, and to allow the reader to compare and contrast with the results in \cite{switch_arxiv}, here we use definitions and notation consistent with \cite{switch_arxiv}.
Since an $n\times n$ switch has $n^2-\text{dimensional}$ queues, we often deal with the Euclidean space $\bR^{n^2}$. However, we represent vectors in $\bR^{n^2}$ as $n\times n$ matrices for convenience. Thus, $\vx$ is a matrix with the $(i,j)^{\text{th}}$ component denoted by $x\lij.$
Thus, for two vectors $\vx$ and $\vy$ in $\bR^{n^2}$, their inner product $\left\langle \vx,\vy\right\rangle $ and Euclidean norm $\|\vx\|$ are defined by
\ba
\left\langle \vx,\vy\right\rangle  \triangleq \sum_{i=1}^n\sum_{j=1}^n x\lij y\lij,
    \quad \|\vx\|\triangleq \sqrt{\left\langle \vx,\vx\right\rangle} = \sqrt{\sum_{i=1}^n\sum_{j=1}^n x\lij^2}.
\ea
For two vectors $\vx$ and $\vy$ in $\bR^{n^2}$, $\vx \leq \vy$ means $x\lij \leq y\lij$ for every $(i,j)$.
We use \vone to denote the all ones vector.
Let \vei denote the vector defined by $\ei\lij =1$ for all $j$ and $\ei_{i'j}=0$ for all $i' \neq i$ and for all $j$. Thus, \vei is a matrix with $i^{\text{th}}$ row being all ones and zeros every where else.
Similarly, let \vetj denote the vector defined by $\etj\lij =1$ for all $i$ and $\etj_{ij'}=0$ for all $j' \neq j$ and for all $i$, i.e., it is a matrix with $j^{\text{th}}$ column being all ones and zeros every where else.
We use $\sum_i(.)$ without the limits on summation to denote $\sum_{i=1}^{n}(.)$.
For a random process $\vq(t)$ and a function $V(.)$, we will sometimes use $V(t)$ to denote $V(\vq(t))$. We use $\text{Var}(.)$ to denote variance of a random variable and $\text{Cov}(.)$ to denote covariance between two random variables.

\subsection{The Switch Model and the MaxWeight Algorithm}\label{sub:model}

For the purposes of queueing-theoretic analysis, an $n\times n$ switch can be treated as an $n\times n$ matrix of queues operating in a time-slotted discrete-time fashion.
Let $\aij(t)$ denote the number of packet arrivals to $(i,j)^{\text{th}}$ queue, i.e., the queue in $i^{\text{th}}$ row and $i^{\text{th}}$ column. We let $\va\in \bR^{n^2}$ denote the vector $(a\lij)\lij$. For every Queue $(i,j)$, the arrival process $\aij(t)$ is a stochastic process
that is i.i.d. across time, with mean $\bE[\aij(t)]=\lambda \lij$ and variance $\text{Var}(\aij(t))=\sigma^2 \lij$ for any time $t$.
We assume that the arrival processes are independent across queues (the processes $a\lij(t)$ and $a_{i'j'}(t)$ are independent for $(i,j)\neq (i',j')$) and are also independent of the queue lengths or schedules chosen in the switch. We further assume that for all $i,j,t$, $\aij(t) \leq a_{\max}$  for some  $a_{\max}\geq 1$.
The arrival rate vector is denoted by $\vlam = (\lambda \lij)\lij$ and the variance vector $(\sigma^2 \lij)\lij$ is denoted by $(\vsig)^2$ or $\vsig^2$. We will use $\vsig$ to denote $(\sigma\lij)\lij$.
We denote the queue length of packets at input port $i$ to be delivered at output port $j$ at time $t$ by $\qijt$. Let $\vq\in \bR^{n^2}$ denote the vector of all queue lengths.

The key scheduling constraints are that (i) at most one packet can be removed from each queue in each time slot and (ii) at most one queue can be served in each row and each column in each time slot. These constraints arise from technological constrains in a real switch, where each row represents an input port and each column represents an output port; see \cite{srikantleibook} for example.
The scheduling constraints can be captured in graph-theoretic language as follows. Let $G$ denote a complete $n \times n $ bipartite graph with $n^2$ edges. Each node on the left side of the bipartite graph can be thought of as representing a row in the matrix of switches and each node on the right side represents a column.
The schedule in each time slot is a matching on this graph $G$. Let $s\lij$ be the amount of service provided to Queue $(i,j)$ in a given time slot. Thus,
$s\lij=1$ if the link between input port $i$ and output port $j$ is matched or scheduled and $s\lij=0$ otherwise and we denote $\vs = (s\lij)\lij$. Then, the set of feasible schedules,
$\cS \subset \{0,1\}^{n^2}$ is the set of all vectors $\vs$ which satisfy
\ba
\sum_{i=1}^{n}s\lij \leq 1, \sum_{j=1}^{n}s\lij \leq 1  \> \forall\> i,j\in\{1,2,\ldots,n\}.
\ea
Let $\cS^*$ denote the set of maximal feasible schedules. Then, it is easy to see that $\cS^*$ is the set of all vectors $\vs$ which satisfy
\ban
\sum_{i=1}^{n}s\lij = 1 , \sum_{j=1}^{n}s\lij = 1   \> \forall\> i,j\in\{1,2,\ldots,n\}
\label{eq:maximal_sched}.
\ean

A scheduling policy or algorithm picks a schedule $\vs(t)$ in every time slot based on the current queue length vector, $\vq(t)$.
In each time slot, the order of events is as follows. Queue lengths at the beginning of time slot $t$ are $\vqt$. A schedule $\vs(t)$ is then picked for that time slot based on the queue lengths. Then, arrivals for that time $\va(t)$ happen. Finally the packets are served and there is unused service if there are no packets in a scheduled queue. The queue lengths are then updated to give the queue lengths for the next time slot. The queue lengths therefore evolve as follows.
\ban
q\lij(t+1)& = \left[\qijt+\aij(t)-\sij(t)\right]^+ \nonumber\\
& = \qijt+\aij(t)-\sij(t)+\uij(t)\nonumber\\
\vq(t+1) & = \vq(t)+\va(t)-\vs(t)+\vu(t) \label{eq:q_evolution},
\ean
where $[x]^+ = \max(0,x)$ is the projection onto positive real axis, $u\lij(t)$ is the unused service on link $(i,j)$. Unused service is $1$ only when link $(i,j)$ is scheduled, but has zero queue length; and it is $0$ in all other cases.
Thus, we have that when $u\lij(t)=1$, we have $\qij(t)=0$, $a\lij(t)=0$, $s\lij(t)=1$ and $\qij(t+1)=0$.
Therefore, we have $u\lij(t)\qijt=0$, $u\lij(t)a\lij(t)=0$ and $u\lij(t)\qij(t+1)=0$.
Also note that since $u\lij(t)\leq s\lij(t)$, we have that $\sum_{i=1}^{n}u\lij \in \{0,1\}$ and $\sum_{j=1}^{n}u\lij\in \{0,1\}$ for all $i,j$.

The MaxWeight Algorithm is a popular scheduling algorithm for switches. In every time slot $t$, each link $(i,j)$ is given a weight equal to its queue length $q\lij(t)$ and the schedule with the maximum weight among the feasible schedules $\mathcal{S}$ is chosen at that time slot.
In other words, using queue lengths as the weights, the permutation matrix with the maximum weight is picked in every time slot.
This algorithm is presented in Algorithm \ref{alg:maxwt}.

\begin{algorithm}
\caption{\label{alg:maxwt}MaxWeight Scheduling Algorithm for an input-queued switch}

Consider the complete bipartite graph described earlier. Let the queue length $q\lij(t)$ be the weight of the edge between input port $i$ and output port $j$. A maximum weight matching in this graph is chosen as the schedule in every time slot, i.e.,
\ban
\vs(t) = \argmax_{\vs \in \mathcal{S}} \sum\lij q\lij(t) s\lij = \argmax_{\vs \in \mathcal{S}} \left \langle \vq(t),\vs  \right\rangle \label{eq:MaxWt}
\ean
Ties are broken uniformly at random.
\end{algorithm}

Note that there is always a maximum weight schedule that is maximal. If the MaxWeight schedule chosen at time $t$, $\vs$ is not maximal, there exists a maximal schedule $\vs^*\in \cS^*$ such that $\vs \leq \vs^*$ . For any link $(i,j)$ such that $s\lij=0$ and $s^*\lij=1$, $\qijt=0$. If not, $\vs$ would not have been a maximum weight schedule. Therefore, we can pretend that the actual schedule chosen is $\vs^*$ and the links $(i,j)$ that are in $\vs^*$ but not in $\vs$ have an unused service of $1$. Note that this does not change the scheduling algorithm, but it is just a notational convenience. Therefore, without loss of generality, we assume that the schedule chosen in each time slot is a maximal schedule, i.e.,
\ban
\vs(t)\in \cS^*   \text{ for all time }t.\nonumber
\ean Hence the MaxWeight schedule picks one of the $n!$ possible permutations from the set $\cS^*$ in each time slot.

Under i.i.d. arrivals, the queue lengths process $\vq(t)$ is a Markov chain. The switch is said to be stable under a scheduling policy if the sum of all the queue lengths is finite in an appropriate stochastic sense (see \cite{srikantleibook} for example).
The capacity region of the switch is the set of arrival rates \vlam for which the switch is stable under some scheduling policy. A policy that stabilizes the switch under any arrival rate in the capacity region is said to be throughput optimal.
It is well known 
\cite{TasEph_92,mckeown96walrand} that the capacity region \cC of the switch is convex hull of all feasible schedules,
\ban
\cC \hn=& \text{Conv}(\mathcal{S}) \nonumber\\
= & \left\{\vlam \in \bR^{n^2}_+ : \sum_{i=1}^{n} \lambda\lij \leq 1  , \sum_{j=1}^{n} \lambda \lij \leq 1   \> \forall\> i,j\in\{1,\ldots,n\} \right\}\nonumber\\
= & \hn\left\{\hn\vlam \hn\in\hn \bR^{n^2}_+ \hn\hn:\hn \La\hn\vlam,\vei \hn\Ra\hn  \leq 1, \La\hn\vlam,\vetj \hn\Ra\hn  \leq 1 \> \forall\> i,j\in\{1,\ldots,n\} \hn\right\}\hn\label{eq:capregion}.\hn\hn
\ean
For any arrival rate vector $\vlam$, $\rho \triangleq \max\lij\{\sum_i\lambda\lij,\sum_j\lambda\lij\}$ is called the load.
It is also known that the queue lengths process is positive recurrent under the MaxWeight algorithm whenever the arrival rate is in the capacity region \cC (equivalently, load $\rho<1$) and therefore is throughput optimal.

For any arrival rate in the capacity region $\cC$, due to positive recurrence of $\vq(t)$, we have that a steady state distribution exists under MaxWeight policy. Let \bvq denote the steady state random vector.
In this paper, we focus on the weighted average queue length under the steady state distribution, i.e., $\bE[\sum_{i,j}\alpha\lij\bqij]$, for some weights $\alpha\lij$, which can be shown to exist as in \cite{switch_arxiv}. We consider a set of switch systems indexed by a parameter $\epsilon$, with arrival rate $\vlam^{\epsilon}$ so that the arrival rate approaches the vector \vnu on the boundary of the capacity region \cC in the limit as $\epsilon \rightarrow 0$. This is called the heavy traffic limit. We are interested in the weighted average queue length in heavy traffic limit, i.e., $\lim_{\epsilon \rightarrow 0} \bE[\sum_{i,j}\alpha\lij\bqij]$. In particular, in this paper, we will consider cases where the sum of the arrival rates at some rows and some columns approach $1,$ and they may approach $1$ at different rates at each column and row.

\subsection{Kingman bound for a discrete-time queue}

To establish our results, we later show that the total queue length along each row and each column is lower bounded. For this purpose, we use a bound on the steady-state queue length in a simple discrete-time queue \cite{erysri-heavytraffic}. While the well-known Kingman bound is for continuous-time G/G/1 queues, due to the similarity in establishing the result, the bound for the discrete-time case is also called the Kingman bound in \cite{erysri-heavytraffic} and we use the same terminology here.  We state the version of the result for the special case where a queue can serve only one packet per slot here since this is what is used in this paper. In this special case, instead of a bound, one has an exact expression for the mean queue length which we present below.
\begin{lemma}\label{lem:Kingman}
  Consider a single server operating in discrete time. In each time slot, packets arrive according to an i.i.d arrival process $\alpha(t)$ with mean $\lambda$ and variance $\sigma^2.$ Let $q$ denote the queue length. Each packet needs exactly one time slot of service. The server operates according to any nonidling policy, \emph{serving one packet in every time slot} whenever the queue is nonempty. Then, the queue is positive recurrent as long as $\lambda<1$, and the steady state mean queue length is given by
  \ba
    \bE[q]=\frac{\sigma^2}{2(1-\lambda)}-\frac{\lambda}{2}.
  \ea
\end{lemma}
We note that the first term on the right-hand side of the above equation is what is referred to as the Kingman bound in \cite{erysri-heavytraffic}.

\subsection{Moment bounds from Lyapunov drift conditions}
In later sections in the paper, we establish state space collapse results by obtaining moment bounds on certain quantities related to the queue length vector based on drift of a Lyapunov function.
A key ingredient in this approach is to obtain moment bounds from drift conditions. A lemma from \cite{hajek_drift} was used in \cite{erysri-heavytraffic} to obtain these bounds and a different result from \cite{bertsimas_momentbound} was used in \cite{switch_arxiv} to obtain tighter bounds. Here we state \cite[Theorem 1]{bertsimas_momentbound} in the form it was stated in \cite{switch_arxiv}.

\begin{lemma} \label{lem:Hajek}
For an irreducible and aperiodic Markov chain $\{X(t)\}_{t\geq 0}$ over a countable state space $\mathcal{X},$ suppose $Z:\mathcal{X}\rightarrow \bR_+$ is a nonnegative-valued Lyapunov function. We define the drift of $Z$ at $X$ as $$\Delta Z(X) \triangleq [Z(X(t+1))-Z(X(t))]\>\mathcal{I}(X(t)=X),$$ where $\mathcal{I}(.)$ is the indicator function. Thus, $\Delta Z(X)$ is a random variable that measures the amount of change in the value of $Z$ in one step, starting from state $X.$ This drift is assumed to satisfy the following conditions:
\begin{enumerate}[label=\textbf{C\arabic*},ref=C.\arabic*] 
\item\label{cond:C1} There exists an $\eta>0,$ and a $\zeta<\infty$ such that for any $t=1,2,\ldots$ and for all $X\in \mathcal{X}$ with $ Z(X)\geq \zeta,$
 \begin{eqnarray*}
    \bE[\Delta Z(X) | X(t)=X]
    \leq - \eta .
 \end{eqnarray*}
\item\label{cond:C2} There exists a $D < \infty$ such that for all $X \in \mathcal{X},$
  \begin{eqnarray*}
  \bP\left(|\Delta Z(X)| \leq D\right) = 1.
  \end{eqnarray*}
\end{enumerate}
Further assume that the Markov chain $\{X(t)\}_t$ converges in distribution to a random variable $\overline{X}$. 
Then, for any $r=1,2,\ldots$,
\ba
\bE[Z\left(\overline{X}\right)^r]\leq &
\left(2\zeta\right)^r+
 \left(4D\right)^r \left(\frac{D+\eta}{\eta} \right)^r r!.
\ea
\end{lemma}

\section{Incompletely Saturated Switch }\label{sec:incomplete}

In this section, we will study the switch system when an arbitrary number of ports are saturated.
We consider the switch where $\no\leq n$ input ports (rows) and $\nt \leq n$ output ports (columns) are saturated.
Without loss of generality, we assume that input ports (rows) $1,2,...,n_1 $  and output ports (columns) $1,2,...,n_2$ are saturated. Thus, we consider a point $\vnu$ on the boundary of the capacity region that lies in $\text{Relint}(\cF \linc)$, the relative interior of the face $\cF \linc$ defined by
\ban
\cF\linc \hn\triangleq&\hn \left\{\hn \vnu \in \cC : \sum_{j=1}^{n} \nu\lij = 1  \> \forall\> i\in\{1,\ldots,n_1\} ,
\sum_{i=1}^{n} \nu \lij = 1\> \forall\> j\in\{1,\ldots,n_2\}\hn\right\} \nonumber\\
= &\hn \left\{\hn\vnu\in \cC \hn:\hn \La\vnu,\vei \Ra  = 1\> \forall\> i\in\{1,\ldots,n_1\} ,\hn
\La\vlam,\vetj \Ra  =1 \> \forall\> j \hn\in\hn \{1,\ldots,n_2\} \hn\right\}\label{eq:face}.
\ean
In other words, if we let $\delta_i =1-\sum_j \nu\lij=1-\La\vnu,\vei \Ra$ and $\tilde{\delta}_j =1-\sum_i \nu\lij=\La\vnu,\vetj \Ra$, we have that  $\delta_i =0$ for $i=1,...n_1$, $\widetilde{\delta}_j =0$ for $j=1,...n_2$ and $\delta_i >0$ for $i>n_1$, $\widetilde{\delta}_j >0$ for $j>n_2$.

We consider a sequence of systems indexed by $\epsilon$. In this section, we consider an i.i.d arrival process $\va\se(t)$ with mean and variance given by
\ba \bE  [\va\se(t)]=& \quad \vlam\se \quad  =\vnu -\epsilon\vk \\
Var[\va\se(t)]=&\quad \left(\vsig^{(\epsilon)}\right)^2 \ea
such that as $\epsilon \rightarrow 0$, $\left(\vsig^{(\epsilon)}\right)^2 \rightarrow \vsig^2$. Here $\vk\in \bR_{+}^{n^2}$ is vector that represents the rates of saturation of different ports. Define
\ban
\kappa_i(\vk)&\triangleq\La\vk,\vei\Ra= \sum_j k\lij  \text{ and } \nonumber\\
\widetilde{\kappa}_j(\vk)&\triangleq\La\vk,\vetj\Ra=\sum_i k\lij.\label{eq:kappas_def}\ean
For simplicity of notation, we will suppress the dependence on $\vk$. Note that $\sum_i \kappa_i=\sum_j \widetilde{\kappa}_j$.
Let $\kappa_{avg} \triangleq \sum_i \kappa_i/n=\sum_j \widetilde{\kappa}_j/n$, $\kappa_{\min}=\min_i\kappa_i$, $\widetilde{\kappa}_{\min}=\min_j \widetilde{\kappa}_j$ and similarly $\kappa_{\max},\widetilde{\kappa}_{\max}$.
Note that in \cite{switch_arxiv}, the setting when $\vk=\vnu$, is studied, in which case  $\kappa_i=1$ and $\widetilde{\kappa_j}= 1$ for all $i,j$. In order to make sure that the heavy traffic parameter $\epsilon$ is comparable to this case, we assume without loss of generality that  $\kappa_{avg}=1$. In other words, we normalize the vector $\vk$ by assuming that $\La \vk,\vone\Ra=n$. 
Such a normalized \vk is called the saturation rate vector.
We will study the switch in the heavy traffic limit as $\epsilon \downarrow 0$.
 Define \ba
\ri\se &\triangleq  1-\sum_j \lambda\se\lij = \delta_i+\epsilon\ki\\
\rtj\se &\triangleq 1-\sum_i \lambda\se\lij = \widetilde{\delta}_j+\epsilon\ktj.
\ea
Note that $\ri\se =\epsilon\ki$ for $i \leq n_1$, and $\rtj\se=\epsilon\ktj $ for $j\leq n_2$. For the unsaturated ports, $\lim_{\epsilon\downarrow0}\ri\se = \delta_i>0$ for $i > n_1$, and $\lim_{\epsilon\downarrow0}\rtj\se=\widetilde{\delta}_j>0 $ for $j>n_2$.

\subsection{Universal Lower Bound}\label{sec:ULB_inc}
We now present  lower bounds on the steady state queue lengths that is satisfied by any scheduling algorithm.
\begin{proposition}\label{prop:Universal_LB}
Consider a set of switch systems with the arrival processes
$\va\se(t)$ described above, parameterized by $0<\epsilon<1$, such that the mean arrival rate vector is $\vlam^{\epsilon}=\vnu-\epsilon \vk$
for some $\vnu \in \text{Relint}(\cF\linc)$,
and the variance is $\left(\vsig^{(\epsilon)}\right)^2$.
For $1\leq i,j \leq n$, \ri, \rtj are defined as above.
Fix a scheduling policy under which the switch system is stable for any $0<\epsilon<1$.
Let $\vq^{(\epsilon)}(t)$ denote the queue lengths process under this policy for each system. Suppose that this process converges in distribution to a steady state random vector $\bvq^{(\epsilon)}$.
Then, for each of these systems, the steady state mean queue lengths can be lower  bounded as follows.
\ban
\bE \left [\sum_j \bqij \se \right ] \geq & \frac{\sum_j \left(\sigma\se_{ij}\right)^2  }{2\ri}- \frac{1-\ri}{2} \text{ for all } 1\leq i \leq n \label{eq:ULBi_incomplete}\\
\bE \left [\sum_i \bqij \se \right ] \geq & \frac{\sum_i \left(\sigma\se_{ij}\right)^2  }{2\rtj}- \frac{1-\rtj}{2} \text{ for all } 1\leq j \leq n \label{eq:ULBj_incomplete}.
\ean
Therefore, in the heavy-traffic limit as $\epsilon \downarrow 0$, if $\left(\vsig^{(\epsilon)}\right)^2 \rightarrow \vsig^2$, for the saturated ports, we have
\ba
\liminf_{\epsilon \downarrow 0} \epsilon \bE \left [\sum_j \bqij \se \right ] \geq
\frac{\sum_j \sigma_{ij}^2  }{2\ki} \text{ for all } 1\leq i \leq n_1\\
\liminf_{\epsilon \downarrow 0} \epsilon \bE \left [\sum_i \bqij \se \right ] \geq
\frac{\sum_i \sigma_{ij}^2  }{2\ktj} \text{ for all } 1\leq i \leq n_2,
\ea
and for the unsaturated ports, we have
\ba
\liminf_{\epsilon \downarrow 0} \epsilon \bE \left [\sum_j \bqij \se \right ] \geq
0 \text{ for all } i > n_1\\
\liminf_{\epsilon \downarrow 0} \epsilon \bE \left [\sum_i \bqij \se \right ] \geq
0 \text{ for all } i > n_2.
\ea
\end{proposition}
\begin{proof}
The queue lengths at each port can be lower bounded by a single server queue as follows. Consider $\sum_j q\lij \se(t)$, the total queue length at input port (row) $i$. It can be lower bounded sample path wise  by a coupled single server queue with arrival process $\sum_j a\lij \se(t)$ \cite[Proposition 1]{switch_arxiv}.
The mean and variance of the arrival process for this single server queue are then $(1-\gamma_i)$ and $\left(\sigma\se_{ij}\right)^2$ respectively because of the independence of the arrival processes across the queues in the matrix. Then, using the Kingman bound for single server queue in Lemma \ref{lem:Kingman},  we get \eqref{eq:ULBi_incomplete}. Similarly lower bounding the total queue length for output port (column) $j$, $\sum_i q\lij \se(t)$ by a single server queue, we get \eqref{eq:ULBj_incomplete}. Taking the heavy traffic limits using the fact that $\ri\se =\epsilon\ki$, $\rtj\se=\epsilon\ktj $ for saturated ports and $\lim_{\epsilon\downarrow0}\ri\se >0$, $\lim_{\epsilon\downarrow0}\rtj\se>0 $ for unsaturated ports, give the heavy traffic limits.
\qed\end{proof}

\subsection{State Space Collapse}\label{sec:SSC_inc}
In this subsection, we will show that under the MaxWeight algorithm, the queue lengths vector concentrates close to a lower dimensional cone. In order to make this more precise, we need to first present the following definitions.

The heavy traffic rate vector \vnu lies in the relative interior of the face $\cF\linc$ which is at the intersection of hyperplanes with the $n_1+n_2$ normal vectors, $\{\vei \text{ for } 1\leq i \leq n_1\}\cup \{\vetj \text{ for } 1\leq j \leq n_2\}$. Call the cone spanned by these normal vector $\cK\linc$, and the subspace spanned by these normal vectors $\cS\linc$ i.e.,
\ba
\cK\linc \triangleq & \Bigl\{ \vx \in \bR^{n^2} \hspace{-.25em}: \vx \hspace{-.25em}=\hspace{-.25em} \sum_{i=1}^{n_1} w_i \vei \hspace{-.25em}+\hspace{-.25em} \sum_{j=1}^{n_2} \widetilde{w}_j \vetj \text{ where }\\
& \quad  w_i \in\bR_+ \for 1\leq i \leq n_1, \widetilde{w}_j \in \bR_+ \text{ for } 1\leq j \leq n_2 \Bigr\} \\
= & \Bigl\{ \vx \in \bR^{n^2} \hspace{-.25em}: \vx \hspace{-.25em}=\hspace{-.25em} \sum_{i=1}^{n} w_i \vei \hspace{-.25em}+\hspace{-.25em} \sum_{j=1}^{n} \widetilde{w}_j \vetj \text{ where }\\
& \quad  w_i \in\bR_+ \for 1\leq i \leq n_1 \text{ and } w_i=0 \for i>n_1,\\
 & \quad \widetilde{w}_j \in \bR_+ \text{ for } 1\leq j \leq n_2 \text{ and } \widetilde{w}_j=0 \for j>n_2\Bigr\}.
\ea
\ban
\cS\linc \triangleq & \Bigl\{ \vx \in \bR^{n^2} \hspace{-.25em}: \vx \hspace{-.25em}=\hspace{-.25em} \sum_{i=1}^{n_1} w_i \vei \hspace{-.25em}+\hspace{-.25em} \sum_{j=1}^{n_2} \widetilde{w}_j \vetj \text{ where } \nonumber\\
& \quad  w_i \in\bR \for 1\leq i \leq n_1, \widetilde{w}_j \in \bR \text{ for } 1\leq j \leq n_2 \Bigr\} \label{eq:space}\\
= & \Bigl\{ \vx \in \bR^{n^2} \hspace{-.25em}: \vx \hspace{-.25em}=\hspace{-.25em} \sum_{i=1}^{n} w_i \vei \hspace{-.25em}+\hspace{-.25em} \sum_{j=1}^{n} \widetilde{w}_j \vetj \text{ where }\nonumber\\
& \quad  w_i \in\bR \for 1\leq i \leq n_1 \text{ and } w_i=0 \for i>n_1,\nonumber\\
& \quad \widetilde{w}_j \in \bR \text{ for } 1\leq j \leq n_2 \text{ and } \widetilde{w}_j=0 \for j>n_2\Bigr\}.\nonumber
\ean
The components of any vector $\vx$ in the subspace $\cS\linc$ can be written in the form, $x\lij=w_i+\wtj$ where $\wi\in \bR$ for $1\leq i \leq n_1$, $\wi=0$ for $ i > n_1 $, $\wtj\in \bR$ for $1\leq j\leq n_2$, and $\wtj=0$ for $ j > n_1 $. The same is true for any vector $\vx$ in the cone $\cK\linc$ with further restriction that $\wi\geq 0 $ for $1\leq i \leq n_1$ and $\wtj\geq0$ for $1\leq j\leq n_2$. This leads to the following lemma relating the structure of the cone $\cK\linc$ and the subspace $\cS\linc$.

\begin{lemma}\label{lem:cone_space}
Let $\no < n$ and $\nt < n$. The cone $\cK\linc$ is the intersection of the space $\cS\linc$ and the positive orthant, i.e.,
\ba
\cK\linc= \cS\linc \cap \bR^{n^2}_+.
\ea
\end{lemma}
\begin{proof}
From the definitions above, it is clear that $\cK\linc \subseteq \cS\linc$ and $\cK\linc \subseteq \bR^{n^2}_+$. Therefore, we have $\cK\linc \subseteq \cS\linc \cap \bR^{n^2}_+$.
Now suppose that $\vx \in \cS\linc \cap \bR^{n^2}_+$, we have that $x\lij = w_i+\widetilde{w}_j \geq 0 $, where $\wi=0$ for $ i > n_1 $ and $\wtj=0$ for $ j > n_1 $. Since $\no < n$, we have $x_{in}=w_i \geq 0$, and so we get that $\wi\geq 0 $ for $1\leq i \leq n_1$. Similarly, we get that $\wtj\geq0$ for $1\leq j\leq n_2$, proving that $\vx \in \cK\linc$ and so $ \cS\linc \cap \bR^{n^2}_+ \subseteq \cK\linc$.
\qed\end{proof}

Let $\vx\parainc$ denote the projection of \vx onto the convex cone $\cK\linc$, and let $\vx\perinc \triangleq \vx-\vx\parainc$ be the perpendicular component. Similarly, let $\vx\parains$ denote the projection of \vx onto the subspace $\cS\linc$, and let $\vx\perins \triangleq \vx-\vx\parains$ be the perpendicular component. For simplicity of notation, we suppress the dependence on \no and \nt.
In \cite{switch_arxiv}, $\vx\para$ and $\vx\per$ were used to denote the projections we denote here by $\vx\parainc$ and $\vx\perinc$ respectively.
We will show that under the MaxWeight Algorithm, all the moments of $\vq\perinc$ are bounded in steady state independent of $\epsilon$. Since the $\ell_1$ norm of the queues length vector, $\|\vq\|_1$ is $\Omega(1/\epsilon)$ as shown in the previous subsection, this establishes that the perpendicular component $\vq\perinc$ is a negligible part of the queue lengths vector $\vq$ for small $\epsilon$. Thus, we establish state space collapse onto the cone $\cK\linc$.

For $\vnu \in \text{Relint}(\cF\linc)$, the vector in the relative interior of the face $\cF\linc$, let $\nu_{min}\triangleq \min_{i,j}\nu_{ij}$. We assume that $\numin > 0$. Then, $\numinp>0$ where
\ba
\numinp \hn\triangleq \min\left\{\hn \nu_{min}, \min_{ i>n_1}\left\{ \hn 1-\La\vnu,\vei\Ra \hn \right\}, \min_{ j>n_2}\left\{\hn 1-\La\vnu,\vetj\Ra \hn \right\}\hn \right\}.
\ea

\begin{proposition} \label{prop:SSC_inc}
Consider a set of switch systems under MaxWeight scheduling algorithm, with the arrival processes
$\va\se(t)$, parameterized by $0<\epsilon<1$ and maximum possible arrivals in any queue $a_{\max}$. The mean arrival rate vector is $\vlam^{\epsilon}=\vnu-\epsilon\vk$ for some $\vnu \in \text{Relint}(\cF\linc)$ such that $\numin\triangleq \min\lij \nu\lij >0$, and a normalized saturation rate vector $\vk\in\bR^{n^2}_+$ such that $\La\vk,\vone\Ra=n$. 
Let the variance $\left(\vsig^{(\epsilon)}\right)^2$ of the arrival process  be such that $\|\vsig^{(\epsilon)}\|^2 \leq \widetilde{\sigma}^2$ for some $\widetilde{\sigma}^2$ not dependent on $\epsilon$. 
Let $\vq^{(\epsilon)}(t)$ denote the queue lengths process of each system, which is positive recurrent. Therefore, the process $\vq^{(\epsilon)}(t)$ converges to a steady state random vector in distribution, which we denote by $\bvq^{(\epsilon)}$.
Then, 
for each system with $0< \epsilon \leq \numinp/2\|\vk\|$, the steady state queue lengths vector satisfies
\ba
\bE \left [\| \bvq\perins \se \|^r\right ] \leq \bE \left [\| \bvq\perinc \se \|^r\right ] \leq (M_r)^r \>\> \forall r \in\{1,2,\ldots\},
\ea
where $\numinp$ is defined as above and $M_r $ is a function of $r,\widetilde{\sigma},\vnu,a_{\max},\numinp$ but independent of $\epsilon$.
\end{proposition}
\begin{proof}
We omit the superscript $\se$ in this proof for simplicity of notation.
For the Markov chain $\vq$, consider the Lyapunov function $\Wperq\triangleq\|\vq\perinc\|$.
We will use Lemma \ref{lem:Hajek} to obtain moment bounds from the drift of $\Wper(.)$. Similar to \cite[Proposition 2]{switch_arxiv}, under the MaxWeight scheduling algorithm, it can be shown that
\ban
\hspace{1em}&\hspace{-1.25em}\bEq{\Delta \Wperq} \nonumber\\
\leq & \frac{1}{2\|\vq\perinc\|}\left(\|\vlam\|^2+\|\vsig\|^2+n -2\epsilon \La\vq\perinc,\vk\Ra \phantom{2\min_{\bf{r}\in \cC}}\right.\nonumber\\
&\left.+ 2\bEq{\La\vq\parainc,\vs(t)-\vnu\Ra}+2\min_{\bf{r}\in \cC} \La\vq,\vnu-\bf{r}\Ra \hn\right) \label{eq:SSCdriftterms}. \hn
\ean
Recall that since $\vs\in\cS^*$, $\La\vei,\vs(t)\Ra=1$ and $\La\vetj,\vs(t)\Ra=1$ for all $i,j,t$. Similarly, since $\vnu \in \cF\linc$, for $1 \leq i\leq n_1$ and $1 \leq j\leq n_2$, we have $\La\vei,\vnu\Ra=1$ and $\La\vetj,\vnu \Ra=1$. By the definition of the cone $\cK\linc$, the vector  $\vq\parainc$ can be written as $\vq\parainc \hspace{-.25em}=\hspace{-.25em} \sum_{i=1}^{n_1} w_i \vei \hspace{-.25em}+\hspace{-.25em} \sum_{j=1}^{n_2} \widetilde{w}_j \vetj$. Putting all these together, we get $\La\vq\parainc,\vs(t)-\vnu\Ra=0$. We now use the following claim to bound the last term in \eqref{eq:SSCdriftterms}.

\begin{claim} \label{claim:SSC}
For any $\vq \in \bR^{n^2}$ and $\vnu \in \text{Relint}(\cF\linc)$ such that $\numin >0$,
\ba
\vnu+\frac{\numinp}{\|\vq\perinc\|}\vq\perinc \in \cC.
\ea
\end{claim}
\begin{proof}
We will verify that $\vnu+\frac{\numinp}{\|\vq\perinc\|}\vq\perinc$ satisfies all the conditions in the definition of \cC in \eqref{eq:capregion}
Note that $\frac{\vq\perinc}{\|\vq\perinc\|}$ is a unit vector along some direction. Since $\nu\lij>\numinp$, clearly, ${\vnu+\frac{\numinp}{\|\vq\perinc\|}\vq\perinc \in \bR^{n^2}_+}$.

It is well known that for any $\vx \in \cK\linc$, $\La \vq\perinc,\vx \Ra \leq 0$. Since
$\vei \in \cK\linc$ for $1\leq i\leq n_1$, we have  $\La \vq\perinc,\vei \Ra \leq 0$. Then, using the fact that $\vnu \in \cF\linc$, we have for $1\leq i\leq n_1$,
\ba
\La \vnu+\frac{\numinp}{||\vq\perinc||}\vq\perinc,\vei \Ra \leq 1.
\ea
For $i>n_1$,
\ba
\La \vnu+\frac{\numinp}{||\vq\perinc||}\vq\perinc,\vei \Ra = & \La \vnu, \vei\Ra+\numinp \La \frac{\vq\perinc}{||\vq\perinc||},\vei\Ra \\
\stackrel{(a)}{\leq} &\La \vnu, \vei\Ra+\nu_{min}' \qquad
\leq \ \ 1 ,
\ea
where (a) follows from Cauchy-Schwartz inequality and the last inequality follows from the definition of $\numinp$. It can similarly be shown that ${\La\vnu+\frac{\numin}{\|\vq\perinc\|}\vq\perinc ,\vetj \Ra \leq 1}$ for $1  \leq j \leq n_2$ as well as $j>n_2$, proving the claim.
\qed\end{proof}
Using the claim, the last term in  \eqref{eq:SSCdriftterms} can be bounded as
\ba
2\min_{\bf{r}\in \cC} \La\vq,\vnu-\bf{r}\Ra \leq & 2\La\vq,\vnu-\left( \vnu+\frac{\numinp}{\|\vq\perinc\|}\vq\perinc \right)\Ra \\
=& -2\La\vq,\frac{\numinp}{\|\vq\perinc\|}\vq\perinc \Ra \\
=& -2 \numinp\|\vq\perinc\|
\ea
Using this in \eqref{eq:SSCdriftterms} and bounding the $-2\epsilon \La\vq\perinc,\vk\Ra $ term using Cauchy-Schwartz inequality, we get
\ba
\hspace{1em}&\hspace{-1.25em}\bEq{\Delta \Wperq} \nonumber\\
\leq  &  \frac{\|\vlam\|^2+\|\vsig\|^2+n}{2\|\vq\perinc\|}-\numinp+\epsilon \|\vk\|\\
\leq  & \frac{\|\vlam\|^2+\|\vsig\|^2+n}{2\|\vq\perinc\|}-\frac{\numinp}{2} \text{ whenever } \epsilon \leq \frac{\numinp}{2\|\vk\|} \\
\leq & - \frac{\numinp}{4} \text{ for all } \vq \text{ such that } \Wperq\geq \frac{2(\|\vlam\|^2+\|\vsig\|^2+n)}{\numinp},
 \ea
Thus condition \ref{cond:C1} is valid with $\zeta = \frac{2(\|\vlam\|^2+\|\vsig\|^2+n)}{\numinp}$ and $\eta = \frac{\numinp}{4}$. Moreover $\zeta$ can be upper bounded by $\zeta \leq \frac{2(\|\vnu\|^2+\|\vsig\|^2+n)}{\numinp}$, an expression that doesn't contain $\epsilon$.
Condition \ref{cond:C2} can be verified using nonexpansivity of projection and the fact that maximum arrivals at every time are $a_{\max}$ \cite{switch_arxiv}.
 Then from Lemma \ref{lem:Hajek}, we get the bound on $\bE \left [\| \bvq\perinc \se \|^r\right ]$ in the proposition. Since $\cK\linc \subseteq \cS\linc$, we  have $\| \bvq\perins \se \|^r \leq \| \bvq\perinc \se \|^r$, completing the proof.
\qed\end{proof}

%

\subsection{Asymptotically tight Upper and Lower bounds under the MaxWeight policy}
In this subsection, we will use the state space collapse result from the previous section to obtain lower and upper bounds on weighted sum of queue lengths  under the MaxWeight algorithm that are tight in heavy traffic limit.
It turns out that the queue length behavior under the MaxWeight algorithm when there is at least one unsaturated port is qualitatively  different from the case when all ports are saturated. The reason for this is discussed in Corollary \ref{cor:complete_sat} in Section \ref{sec:discussion}. So, in this subsection, we focus on the case when at least one port in not saturated.
If all the input ports are saturated, from the definition of the capacity region, it follows that all the output ports are also saturated, i.e., whenever $\no=n$, we also have $\nt=n$. Similarly, if all the output ports are saturated, it again follows that all the input ports are saturated.
Since we are interested in incompletely saturated switch, in this section, we assume $\no< n$ and $\nt< n$.

The queue length bounds are obtained by setting the drift of the following function to zero in steady state.
\ba
\Vq=\norm{\vq\paras}^2.
\ea
Its drift is defined as
\ba
\Delta \Vq \triangleq & [V(\vq(t+1))-V(\vq(t))]\>\mathcal{I}(\vq(t)=\vq).
\ea
We now state the main result of the paper in a general form. In Section \ref{sec:discussion}, we will interpret this result as well as present various special cases.

\begin{theorem}\label{thm:UBM_inc}
Consider a set of switch systems under MaxWeight scheduling algorithm, with the arrival processes
$\va\se(t)$, parameterized by $0<\epsilon<1$ and maximum possible arrivals in any queue $a_{\max}$. The mean arrival rate vector is $\vlam^{\epsilon}=\vnu-\epsilon\vk$ for some $\vnu \in \text{Relint}(\cF\linc)$ such that $\numin\triangleq \min\lij \nu\lij >0$, and a normalized saturation rate vector $\vk\in\bR^{n^2}_+$ such that $\La\vk,\vone\Ra=n$.
Let the variance $\left(\vsig^{(\epsilon)}\right)^2$ of the arrival process  be such that $\|\vsig^{(\epsilon)}\|^2 \leq \widetilde{\sigma}^2$ for some $\widetilde{\sigma}^2$ not dependent on $\epsilon$ and assume that $\no<n$ and $\nt<n$.
Let $\vq^{(\epsilon)}(t)$ denote the queue lengths process of each system, which is positive recurrent. Therefore, the process $\vq^{(\epsilon)}(t)$ converges to a steady state random vector in distribution, which we denote by $\bvq^{(\epsilon)}$.
Then, 
for each system with $0< \epsilon \leq \numinp/2\|\vk\|$, the steady state queue lengths vector satisfies
\begin{alignat*}{3}
\frac{1}{2\epsilon}\La \left(\vsig\se\right)^2, \vzeta  \Ra-B_1(\epsilon)
 &\leq\bE\left[\La\bvq\se,\valpha \Ra\right]
 &&\leq  \frac{1}{2\epsilon}\La \left(\vsig\se\right)^2, \vzeta  \Ra+B_2(\epsilon)
\end{alignat*}
for any fixed weight vector $\valpha \in \bR^{n^2}$ such that $\La\valpha,\vei \Ra=n\ki$ for $ i \leq \no$ and $\La\valpha,\vetj \Ra=n\ktj$ for $j \leq \nt$, where  $B_1(\epsilon)$ as well as $B_2(\epsilon)$ are $o(\frac{1}{\epsilon})$, i.e.,  $\lim_{\epsilon\rightarrow 0} \epsilon B_1(\epsilon)=0$ and $\lim_{\epsilon\rightarrow 0} \epsilon B_2(\epsilon)=0$. The vector \vzeta is defined by
\ban
\zeta\lij & \triangleq   \left\{
  \begin{array}{@{}ll@{}}
    2-\frac{2n-\no-\nt}{n^2-\no\nt} & \text{if}\ i\leq \no \text{ and } j\leq \nt \\
    1+\frac{\nt}{n^2-\no\nt} & \text{if}\ i \leq \no \text{ and } j >\nt \\
    1+\frac{\no}{n^2-\no\nt} & \text{if}\ i >\no \text{ and } j \leq \nt \\
    0 & \text{if}\ i> \no \text{ and } j> \nt.
  \end{array}\right.\label{eq:zeta_def}
\ean
Thus, in the heavy traffic limit as $\epsilon\downarrow 0$, we have
\ba
\lim_{\epsilon\to 0} \epsilon \bE\left[\La\bvq\se,\valpha \Ra\right] &= \frac{1}{2}\La \vsig^2, \vzeta  \Ra.
\ea
Moreover, for any $i>n_1$ and $j>n_2$,
\ba
\lim_{\epsilon \rightarrow 0} \epsilon\bE\left[\bqij\se\right]=0.
\ea
\end{theorem}
Note that in general, the weights $\alpha\lij$ are allowed to be negative. We now present the proof of the theorem.
\begin{proof}
We consider the switch for a fixed $0< \epsilon \leq \numinp/2\|\vk\|$. For simplicity of notation,  we again omit the superscript $(\epsilon)$ in this proof. Similar to the notation in \cite{switch_arxiv}, we use $\bvq$ to denote the steady state queue length vector  and $\bva$ to denote the steady state arrival vector which is identically distributed to the  vector $\va(t)$ at any time $t$.
We use $\vs(\bvq)$ and $\vu(\bvq)$ for the schedule and unused service to explicitly show their dependence on the queue lengths. 
If the queue length at time $t$ is $\bvq$, then the queue length at time $t+1$, $\bvq+\bva-\vs(\bvq)+\vu(\bvq)$ is denoted by
$\bvq^+$. Since $\bvq$ is the steady state queue length, it has the same distribution as $\bvq^+$.

It can be easily shown using Lemma \ref{lem:Hajek} that in steady state, $\bE[\norm{\bvq}^2]$, is finite and consequently we have
\ban \bE[V(\bvq)] < \infty \quad \text{ and } \quad
\bE[\norm{\bvq}_1] = \bE[\sum\lij\bqij] < \infty, \label{eq:V_finite}.
\ean
where $\norm{.}_1$ denotes the $\ell_1$ norm.
See Lemma 5 in \cite{switch_arxiv} for details.
Setting the drift of $V(\bvq)$ to zero in steady state, we get
\ba
0& = \lefteqn{\bE[\Delta \V(\bvq)]}\\
&= \bE[\V(\bvq+\bva-\vs(\bvq)+\vu(\bvq))-\Vi(\bvq)]\\
&= \bE\left[\norm{\left(\bvq+\bva-\vs(\bvq)+\vu(\bvq)\right)\paras}^2-\norm{\bvq\paras}^2\right] \\
&\stackrel{(a)}{=} \bE\left[\norm{\bvq\paras+(\bva-\vs(\bvq))\paras+\vu\paras(\bvq)}^2-\norm{\bvq\paras}^2\right] \\
&= \bE\left[\norm{\bvq\paras+(\bva-\vs(\bvq))\paras}^2
+2\La\bvq\paras+(\bva-\vs(\bvq))\paras, \vu\paras(\bvq)\Ra\right]\\
&\qquad+\bE\left[\norm{\vu\paras(\bvq)}^2-\norm{\bvq\paras}^2\right]\\
&= \bE\left[\norm{\bvq\paras+(\bva-\vs(\bvq))\paras}^2
+2\La\bvq\paras+(\bva-\vs(\bvq))\paras+\vu\paras(\bvq), \vu\paras(\bvq)\Ra\right]\\
&\qquad-\bE\left[\norm{\vu\paras(\bvq)}^2+\norm{\bvq\paras}^2\right]\\
&\stackrel{(b)}{=} \bE\left[\norm{\bvq\paras+(\bva-\vs(\bvq))\paras}^2
+2\La\bvq\paras^+, \vu\paras(\bvq)\Ra-\norm{\vu\paras(\bvq)}^2-\norm{\bvq\paras}^2\right]\\
&= \bE\left[\norm{(\bva-\vs(\bvq))\paras}^2+2\La \bvq\paras, (\bva-\vs(\bvq))\paras\Ra-\norm{\vu\paras(\bvq)}^2+2\La\bvq\paras^+, \vu\paras(\bvq)\Ra\right]
\ea
where (a) and (b) follow from the fact that projection onto a subspace is linear. Therefore, we have
\ban
2\bE\left[\La \bvq\paras, (\vs(\bvq)-\bva)\paras\Ra\right]& = \bE\left[\norm{(\bva-\vs(\bvq))\paras}^2\right]-\bE\left[\norm{\vu\paras(\bvq)}^2\right]\label{eq:terms1_3}\\
&\qquad +2\bE\left[\La\bvq\paras^+, \vu\paras(\bvq)\Ra\right]. \label{eq:term4}
\ean
We will now study each of the terms in this equation. Consider the LHS term in \eqref{eq:terms1_3}. Since any vector of the form $\vx\perins$ is orthogonal to the space $\cS\linc$, we get
\ba
&2\bE\left[\La \bvq\paras, (\vs(\bvq)-\bva)\paras\Ra\right]\\
&=2\bE\left[\La \bvq\paras, (\vs(\bvq)-\bva)\paras\Ra\right]+2\bE\left[\La \bvq\paras, (\vs(\bvq)-\bva)\perins\Ra\right]\\
&= 2\bE\left[\La \bvq\paras, \vs(\bvq)-\bva\Ra\right]\\
&\stackrel{(a)}{=} 2\bE\left[\La \bvq\paras, \vs(\bvq)-\vlam\Ra\right]\\
&= 2\bE\left[\La \bvq\paras, \vs(\bvq)-(\vnu -\epsilon\vk)\Ra\right] \\
&= 2\epsilon \bE\left[\La \bvq\paras, \vk \Ra \right]+2\bE\left[\La \bvq\paras, \vs (\bvq)-\vnu \Ra\right]\\
&\stackrel{(b)}{=} 2\epsilon \bE\left[\La \sum_{i=1}^{n_1} w_i \vei + \sum_{j=1}^{n_2} \widetilde{w}_j \vetj , \vk \Ra \right]  +2\bE\left[\La \sum_{i=1}^{n_1} w_i \vei + \sum_{j=1}^{n_2} \widetilde{w}_j \vetj , \vs (\bvq)-\vnu \Ra\right]\\
&\stackrel{(c)}{=} 2\epsilon \bE\left[\sum_{i=1}^{n_1} w_i\La  \vei,\vk\Ra + \sum_{j=1}^{n_2} \widetilde{w}_j \La \vetj , \vk \Ra \right]\\
&\qquad +2\bE\left[\sum_{i=1}^{n_1} w_i \La \vei ,\vs (\bvq)-\vnu \Ra + \sum_{j=1}^{n_2} \widetilde{w}_j \La \vetj , \vs (\bvq)-\vnu \Ra\right]\\
&\stackrel{(d)}{=} \frac{2\epsilon}{n} \bE\left[\sum_{i=1}^{n_1} w_i\La  \vei,\valpha\Ra + \sum_{j=1}^{n_2} \widetilde{w}_j \La \vetj , \valpha \Ra \right]\\
&= \frac{2\epsilon}{n} \bE\left[\La \sum_{i=1}^{n_1} w_i  \vei +\sum_{j=1}^{n_2} \widetilde{w}_j \vetj , \valpha \Ra \right]\\
&= \frac{2\epsilon}{n} \bE\left[\La \bvq\paras , \valpha \Ra \right]\\
&= \frac{2\epsilon}{n} \bE\left[\La \bvq , \valpha \Ra \right]-\frac{2\epsilon}{n} \bE\left[\La \bvq\perins , \valpha \Ra \right], \numberthis \label{eq:Ti_qperp}
\ea
where (a) follows from the fact that the arrivals are independent of queue lengths. From the definition of the space $\cS\linc$ in \eqref{eq:space}, we know that the vector $\bvq\paras$ can be represented as $\sum_{i=1}^{n_1} w_i \vei + \sum_{j=1}^{n_2} \widetilde{w}_j \vetj$ for some $w_i\in \bR$, $1\leq i \leq \no$ and $\wtj \in \bR$, $1\leq i \leq \no$, giving (b).
Since the schedule is always assumed to be maximal, from \eqref{eq:maximal_sched}, we have that $\La \vei ,\vs\Ra=1$ for all $i$. Since the first \no rows are saturated, we have from \eqref{eq:face} that $\La \vei ,\vnu\Ra=1$ for $1\leq i\leq \no$. Therefore, we get that $\La \vei ,\vs (\bvq)-\vnu \Ra=0$ for $1\leq i\leq \no$, and similarly, we have $\La \vetj , \vs (\bvq)-\vnu \Ra=0$ for $1\leq j\leq \nt$. Consequently, the second term in (c) vanishes.
From the definition of $\kappa_i$ in \eqref{eq:kappas_def} and the assumption on $\valpha$, we have that $\La  \vei,\vk\Ra=\kappa_i=\La  \vei,\valpha\Ra/n$ for $1\leq i\leq \no$. Similarly for $1\leq j\leq \nt$, we have  $\La  \vetj,\vk\Ra=\kappa_i=\La  \vetj,\valpha\Ra/n$ giving us (d). Using Cauchy-Schwartz inequality, we can bound the last term in \eqref{eq:Ti_qperp} as follows,
\begin{alignat}{3}
 - \bE\left[\norm{ \bvq\perins} \right]\norm{\valpha}
&\leq  \bE\left[\La \bvq\perins , \valpha \Ra \right]
&&\leq  \bE\left[\norm{ \bvq\perins} \right]\norm{\valpha} \nonumber
\\%
 -M_1 \norm{\valpha}
&\leq  \bE\left[\La \bvq\perins , \valpha \Ra \right]
&&\leq   M_1\norm{\valpha}, \label{eq:useSSC1}
\end{alignat}
where the last set of inequalities follow from the state space collapse in Proposition \ref{prop:SSC_inc}. Putting this back in \eqref{eq:Ti_qperp}, we get
\begin{alignat}{3}
-\frac{2\epsilon}{n} M_1\norm{\valpha}
&\leq2\bE\left[\La \bvq\paras, (\vs(\bvq)-\bva)\paras\Ra\right]-\frac{2\epsilon}{n} \bE\left[\La \bvq , \valpha \Ra \right]
&&\leq \frac{2\epsilon}{n} M_1\norm{\valpha} \label{eq:Ti}
\end{alignat}

We now consider the first term on the RHS of \eqref{eq:terms1_3}.
Let $\vf_1, \vf_2, ..., \vf_L$ be an orthonormal basis of the space $\cS\linc$ where $L$ is the dimension of the space $\cS\linc$. From the definition of the space $\cS\linc$, we know that each of these vectors $\vf_l$ can be written as
$\flij= \nuli+\nutlj$ for some $\nuli, \nutlj\in \bR$ for all $i,j$ with $\nuli=0$ for $i>n_1$, $\nutlj=0$ for $j>n_2$. Then the norm of projection onto the subspace $\cS\linc$ can be written in terms of the projections onto the basis vectors as follows.
\ba
& \bE\left[\norm{(\bva-\vs(\bvq))\paras}^2\right]\\
&=  \bE\left[\sum_l\La \bva-\vs(\bvq), \vf_l \Ra^2\right] \\
&= \sum_l \bE\left[\left(\sum_{ij}(\baij-\sij(\bvq))f\lij \right)^2\right] \\
&= \sum_l \bE\left[\left(\sum_{ij}(\baij-\sij(\bvq))(\nuli+\nutlj) \right)^2\right] \\
&= \sum_l \bE\left[\left(\sum_{i}\nuli\left(\sum_j (\baij-\sij(\bvq))\right)+\sum_{j}\nutlj\left(\sum_i(\baij-\sij(\bvq))\right) \right)^2\right] \\
&\stackrel{(a)}{=} \sum_l \bE\left[\left(\sum_{i}\nuli\left(\sum_j \baij-(1-\ki\epsilon)\right) +\sum_{j}\nutlj\left(\sum_i\baij-(1-\ktj\epsilon)\right)\right.\right. \\
&\qquad \left.\left.-\epsilon \left(\sum_i \ki\nuli+\sum_j\ktj\nutlj\right)\right)^2\right] \\
&\stackrel{(b)}{=} \sum_lVar\left(\sum_i \nuli\sum_j \baij+\sum_j \nutlj\sum_i \baij\right)+\epsilon^2\sum_l\left(\sum_i \ki\nuli+\sum_j\ktj\nutlj\right)^2\\
&\stackrel{(c)}{=} \sum_l\left[Var\left(\sum_i \nuli\sum_j \baij\right)+ Var\left(\sum_j \nutlj\sum_i \baij\right)\right.\\
&\quad +\left.2Cov\left(\sum_i \nuli\sum_{j'} \overline{a}_{ij'}, \sum_j \nutlj\sum_{i'} \overline{a}_{i'j} \right)\right]+\epsilon^2\sum_l\La\vk,\vf_l\Ra^2\\
&\stackrel{(d)}{=} \sum_l\left[\sum_i \nuli^2 Var\left(\sum_j \baij \right)+\sum_j \nutlj^2 Var\left(\sum_i \baij \right) \right.\\
&\quad +\left.2\sum_{ij} \nuli\nutlj Cov\left(\sum_{j'}\overline{a}_{ij'}, \sum_{i'} \overline{a}_{i'j}\right)\right]+\epsilon^2\norm{\vk\paras}^2\\
&\stackrel{(e)}{=} \sum_l\left[\sum_i \nuli^2 \sum_j \sigma\lij^2+\sum_j \nutlj^2 \sum_i \sigma\lij^2+2\sum_{ij}\nuli\nutlj\sigma\lij^2\right]+\epsilon^2\norm{\vk\paras}^2\\
&= \sum_l\left[\sum_{ij} (\nuli+\nutlj)^2\sigma\lij^2 \right]+\epsilon^2\norm{\vk\paras}^2\\
&= \sum_{ij}\sum_l \flij^2\sigma\lij^2+\epsilon^2\norm{\vk\paras}^2\\
&= \sum_{ij}\sum_l \La \vchi^{(ij)}, \vf_l\Ra \sigma\lij^2+\epsilon^2\norm{\vk\paras}^2\\
&= \sum_{ij}\norm{ \vchi^{(ij)}\paras}^2 \sigma\lij^2+\epsilon^2\norm{\vk\paras}^2\\
& = \frac{1}{n}\La \vsig^2, \vzeta  \Ra +\epsilon^2\norm{\vk\paras}^2 \numberthis\label{eq:Tii},
\ea
where $\vchi^{(ij)}$ is the matrix with $1$ in $(i,j)^{\text{th}}$ position and $0$ everywhere else.
Since a maximal schedule is always picked, from \eqref{eq:maximal_sched}, we get (a). Equation b) follows from the fact that the total arrival rate for row $i$ is $(1-\ki\epsilon)$ and that for column $j$ is $(1-\ktj\epsilon)$. Independence of the arrival process across the ports along with the definition of $\kappa_i$ and $\widetilde{\kappa}_j$ in \eqref{eq:kappas_def} gives (c). We again use the independence of arrival processes to get (d) and (e). The following lemma, which is proved in Appendix \ref{appendix:projections} gives \eqref{eq:Tii}.
\begin{lemma}\label{lem:projections}
For all $1\leq i,j \leq n$,
\ba
\norm{\vchi^{(ij)}\paras}^2 &= \frac{\zeta\lij}{n},
\ea
where $\vzeta$ is defined in \eqref{eq:zeta_def}.
\end{lemma}
%
%
%

We will now focus on the second term on RHS of \eqref{eq:terms1_3}. In order to do this, we will first obtain a  bound on the unused service. We know from \eqref{eq:V_finite} that $\bE[\sum\lij\bqij]$ is finite and so $\bE[\sum_{j}\bqij]$ finite for each $i$. For $i\leq \no$, setting the drift of $\sum_{j}\bqij$ to zero in steady state, we get
\ba
\bE\left[\sum_j \bqij\right] &= \bE\left[\sum_j \bqij^+\right]\\
&= \bE\left[\sum_j\bqij+ \baij -\sij(\bvq)+\uij(\bvq) \right]\\
0&= \sum_j (\nuij -\epsilon k\lij)-1+\bE\left[\sum_j \uij(\bvq)\right]\\
\bE\left[\sum_j\uij(\bvq)\right]&= \epsilon\ki,
\ea
where the last equality follows from \eqref{eq:face} since the row $i$ is saturating. Similarly, for $j\leq n_2$, we have that $\bE\left[\sum_i \uij(\bvq)\right]=\epsilon\ktj$.

For any $\vx \in \bR^{n^2}$, let $\widehat{\vx}\in \bR^{n^2}$ denote its projection on to the space spanned by the vectors $\vchi^{(ij)}$ for $i\leq \no$ or $j \leq \nt$. Call this space $\mathcal{X}_{\no\nt}$. Clearly, this space contains the space $\cS\linc$, i.e., $\cS\linc\subseteq\mathcal{X}_{\no\nt}$. In other words, the vector $\widehat{\vx}\in \bR^{n^2}$ is obtained by replacing the $(n-n_1)(n-n_2)$ components with $i>n_1$ and $j>n_2$ with zeros, i.e.,
\ba
\widehat{x}\lij \triangleq\left\{
  \begin{array}{@{}ll@{}}
    x\lij & \text{if}\ i\leq n_1 \text{ or } j\leq n_2\\
        0 & \text{if}\ i> n_1 \text{ and } j> n_2\\
  \end{array}\right.
\ea
Moreover, for any $\vx\in \bR^{n^2}$, $\vx-\widehat{\vx}$ is orthogonal to the space $\mathcal{X}_{\no\nt}$ and so, is also orthogonal to the space $\cS\linc$.

Now the second term on RHS of \eqref{eq:terms1_3} can be upper bounded as follows.
\ba
\bE\left[\norm{\vu\paras(\bvq)}^2\right]  &= \bE\left[\norm{\left(\widehat{\vu}(\bvq)+\vu(\bvq)-\widehat{\vu}(\bvq)\right)\paras}^2\right] \\
 & \stackrel{(a)}{=} \bE\left[\norm{\left(\widehat{\vu}(\bvq)\right)\paras+\left(\vu(\bvq)-\widehat{\vu}(\bvq)\right)\paras}^2\right] \\
& \stackrel{(b)}{=} \bE\left[\norm{\left(\widehat{\vu}(\bvq)\right)\paras}^2\right] \\
 & \stackrel{(c)}{\leq} \bE\left[\norm{\left(\widehat{\vu}(\bvq)\right)}^2\right] \\
&= \bE\left[\sum\lij\widehat{u}^2\lij(\bvq)\right]\\
& \stackrel{(d)}{=} \bE\left[\sum\lij\widehat{u}\lij(\bvq)\right]\\
&\leq \bE\left[\sum_{i=1}^{n_1}\sum_j \uij(\bvq)\right]+\bE\left[\sum_{j=1}^{n_2}\sum_i \uij(\bvq)\right]\\
& \stackrel{(e)}{=} \epsilon\sum_{i=1}^{n_1}\ki+\epsilon\sum_{j=1}^{n_2} \ktj\\
&\leq 2\epsilon \La\vk,\vone\Ra\\
&=2\epsilon n, \numberthis \label{eq:u_squared}
\ea
where (a) follows from linearity of projection onto a subspace.
Since for any vector $\vx$, $\vx-\widehat{\vx}$ is orthogonal to the space $\cS\linc$, we get (b).
Inequality (c) is true due to the nonexpansive property of projection. Since $\uij\in\{0,1\}$, we have (d). Since the saturation rate vector \vk is assumed to be normalized, we get \eqref{eq:u_squared}. Using the trivial lower bound of zero, we have
\ban
0\leq\bE\left[\norm{\vu\paras(\bvq)}^2\right]  &\leq 2\epsilon n \label{eq:Tiii}.
\ean
We will now consider the final term, the one in \eqref{eq:term4}.
\ba
2\bE\left[\La\bvq\paras^+, \vu\paras(\bvq)\Ra\right]&=2\bE\left[\La\bvq\paras^+, \vu(\bvq)\Ra\right]\\
&\stackrel{(a)}{=}2\bE\left[\La\bvq\paras^+, \widehat{\vu}(\bvq)\Ra\right]\\
&\stackrel{(b)}{=}2\bE\left[\La\bvq^+, \widehat{\vu}(\bvq)\Ra\right]- 2\bE\left[\La\bvq\perins^+, \widehat{\vu}(\bvq)\Ra\right]\\
&= -2\bE\left[\La\bvq\perins^+, \widehat{\vu}(\bvq)\Ra\right], \numberthis \label{eq:quterm}
\ea
where (a) follows from the fact that  $\vx-\widehat{\vx}$ is orthogonal to the space $\cS\linc$. Since by the definition of unused service in \eqref{eq:q_evolution}, we have that $\bqij^+=0$ whenever $\uij(\bvq)>0$, and so the first term in (b) vanishes. Therefore, using Cauchy-Schwartz inequality, we get
\begin{alignat}{3}
-2\sqrt{\bE\left[\norm{\bvq\perins^+}^2\right]\bE\left[\norm{\widehat{\vu}(\bvq)}^2\right]}
&\leq 2\bE\left[\La\bvq\paras^+, \vu\paras(\bvq)\Ra\right]
&& \leq 2\sqrt{\bE\left[\norm{\bvq\perins^+}^2\right]\bE\left[\norm{\widehat{\vu}(\bvq)}^2\right]} \nonumber\\
-2M_2\sqrt{\bE\left[\norm{\widehat{\vu}}^2\right]}
&\leq 2\bE\left[\La\bvq\paras^+, \vu\paras(\bvq)\Ra\right]
&&\leq 2M_2\sqrt{\bE\left[\norm{\widehat{\vu}}^2\right]} \label{eq:useSSC2}\\
-2M_2\sqrt{2\epsilon n}
&\leq 2\bE\left[\La\bvq\paras^+, \vu\paras(\bvq)\Ra\right]
&&\leq 2M_2\sqrt{2\epsilon n} \label{eq:Tiv}
\end{alignat}
where \eqref{eq:useSSC2} is obtained by  using the fact that in steady state, $\bE\left[\norm{\bvq\perins^+}^2\right]=\bE\left[\norm{\bvq\perins}^2\right]$, which is bounded by $M_2$ from state space collapse in Proposition \ref{prop:SSC_inc}. We get \eqref{eq:Tiv} from the bound in \eqref{eq:u_squared}.

Substituting \eqref{eq:Ti}, \eqref{eq:Tii}, \eqref{eq:Tiii} and \eqref{eq:Tiv} in \eqref{eq:terms1_3} and \eqref{eq:term4}, we get the theorem with
\ba
B_1(\epsilon)=  & M_1\norm{\valpha}-\frac{n\epsilon}{2}\norm{\vk\paras}^2+n^2+2M_2\frac{n\sqrt{n}}{\sqrt{2\epsilon}}\\
B_2(\epsilon)=  & M_1\norm{\valpha}+\frac{n\epsilon}{2}\norm{\vk\paras}^2+2M_2\frac{n\sqrt{n}}{\sqrt{2\epsilon}}.
\ea
\qed\end{proof}

The main idea of the proof is to set the drift of a carefully chosen test function $V(.)$ to zero. The choice of this function is crucial to obtain tight heavy traffic bounds. We will now briefly motivate our choice of the function $V(.)$. For a discrete time single server (G/G/1) queue, $q(t)$ that evolves according to $q(t+1)=q(t)+a(t)-s(t)+u(t)$, the right test function to obtain tight queue length bounds is $q^2$. Such a bound is known as Kingman bound \cite[Section 10.1]{srikantleibook}.
Next, consider a load balancing system under the `Join the shortest queue'(JSQ) policy operating in discrete time. There are a finite number of servers, each with a separate queue, similar to super market checkout lanes. Whenever a user arrives into the system, (s)he joins the queue with the shortest length, breaking ties uniformly at random.
Tight heavy traffic queue length bounds are obtained for this system in \cite{erysri-heavytraffic} by first showing that the queue lengths collapse to a single dimension where they are all equal.
Then, tight bounds are obtained by setting the drift of the quadratic test function $\bE[(\sum_i \bq_i)^2]$ to zero in steady state. This function is same as 
$\norm{\vq_{\parallel}}^2$ (up to a factor of $n$, which is not important) where $\vq_{\parallel}$ denotes the projection of the queue length vector \vq onto the region of state space collapse, which is the line along the vector that has one's in all components.

These examples, motivate us to choose the norm square of projection of queue lengths vector onto the region of state space collapse as the test function in general. In the case of the switch, it would be $\norm{\vq_{\parallel \cK}}^2$. However, projection operator onto a convex cone is difficult to study because it is not linear. Moreover, closed form expression for projection onto a general cone is not known. Therefore, we relax this function and use $\norm{\vq\paras}^2$ as the test function.
Since we use the relaxed function, it is sufficient to use state space collapse into the space $\cS\linc$ as opposed to the stronger form of state space collapse into the cone $\cK\linc$.
Note that in the proof of Theorem \ref{thm:UBM_inc}, we use state space collapse only at two instances, viz., \eqref{eq:useSSC1} and \eqref{eq:useSSC2}, and both these use the weaker result, $\bE \left [\| \bvq\perins \se \|^r\right ]$. In other words, we only use the fact that the state collapses to the space $\cS\linc$ and not the cone $\cK\linc$.
Such a relaxation works because of the property of the cone $\cK\linc$ that, it is just the intersection of the space it spans $\cS\linc$ with the positive orthant $\bR^{n^2}_+$ as proved in Lemma \ref{lem:cone_space}. Therefore, any positive queue length vector that is in the space $\cS\linc$ is also in the cone $\cK\linc$. Since queue length vectors are nonnegative, if we know that a queue length vector collapses onto the space $\cS\linc$, we know that it should collapse onto the cone $\cK\linc$.

Even though we only need the weaker version of state space collapse, viz., the bound on $\bE \left [\| \bvq\perins \se \|^r\right ]$, we proved a stronger version, viz., the bound on $\bE \left [\| \bvq\perinc \se \|^r\right ]$ in Proposition \ref{prop:SSC_inc} for completeness.
The proof of weaker state space collapse can be much simpler  because projection onto a subspace is a linear operator while projection onto a convex cone is not. Therefore, we can write down the complete proof of Theorem \ref{thm:UBM_inc} (including the weaker version of Proposition \ref{prop:SSC_inc}) without even referring to the cone $\cK\linc$.

A major contribution of this paper is that we present the test function in a form so that it may be easily generalized to other problem settings beyond just the switch system.
For any problem that exhibits state space collapse into a cone that satisfies a property similar to Lemma \ref{lem:cone_space}, we may be able to use $\norm{\vq\paras}^2$ as the test function. Moreover, we don't need to know the exact closed form $\vq\paras$ to use this test function.
This is because of the following reason.
In general, we expect that the right test function has the form $\norm{\vq\para}^2$ for appropriately defined projection onto a region into which the state collapses. This is because of the following reason.
It was argued in \cite{switch_arxiv} that a quadratic test function should be chosen so that, when its drift is set to zero, among the various terms, the cross term between $\bvq^+$ and $\vu(\bvq)$ should be small under state space collapse. When $\norm{\vq\para}^2$ is used as the test function, such cross terms will be of the form  of the expression in \eqref{eq:term4}. Then these cross terms can be shown to be small under state space collapse using an argument similar to that in the set of equations ending in \eqref{eq:quterm}.

We now make a short remark about the scheduling policy. In the proof of Theorem \ref{thm:UBM_inc}, we do not use any details about the scheduling policy, except for the fact that it is throughput optimal and exhibits state space collapse as in Proposition \ref{prop:SSC_inc}. Therefore, the queue length bounds in Theorem \ref{thm:UBM_inc} hold true under any throughput optimal algorithm that exhibits state space collapse as in Proposition \ref{prop:SSC_inc}. Moreover, in the proof of
state space collapse   in Proposition \ref{prop:SSC_inc}, we use the details of MaxWeight scheduling policy only in 
the Claim \ref{claim:SSC}.

\section{Discussion}\label{sec:discussion}

In this section, we present various corollaries and extensions of Theorem \ref{thm:UBM_inc} and interpret the results.
The following corollary gives a bound on sum of the queue lengths.
\begin{corollary}\label{cor:kappa_min}
Consider the set of switch systems operating under the MaxWeight algorithm  as described in Theorem \ref{thm:UBM_inc} with $0\leq \no,\nt<n$. Then, in the heavy traffic limit, we have
\ba
\frac{1}{2\kappa_{\max}}\La \vsig^2, \vzeta  \Ra\leq \lim_{\epsilon\to 0} \epsilon \bE\left[\sum\lij\bq\lij\se\right] &\leq  \frac{1}{2\kappa_{\min}}\La \vsig^2, \vzeta  \Ra,
\ea
where $\kappa_{\max} =\max_{i\leq\no,j\leq\nt}\{\ki,\ktj\}$ and  $\kappa_{i\leq\no,j\leq\nt} =\min\lij\{\ki,\ktj\}$. Moreover, under any stable algorithm the sum queue lengths is lower bounded by
\ba
\lim_{\epsilon\to 0} \epsilon \bE\left[\sum_{ij} \bqij\se\right] &\geq  \frac{1}{2}\max\left\{\La \vsig^2, \vzeta'  \Ra,\La \vsig^2, \vzeta''  \Ra\right\},
\ea
where $\zeta'\lij = 1/\ki$ and $\zeta'\lij = 1/\ktj$.
\end{corollary}
\begin{proof}
It is easy to see that there exists a weight vector $\valpha \in \bR^{n^2}_+$  such that $\kappa_{\min}\leq\alpha\lij$ for all $i,j$, $\La\valpha,\vei \Ra=n\ki$ for $ i \leq \no$ and $\La\valpha,\vetj \Ra=n\ktj$ for $j \leq \nt$. Using such a weight vector in  Theorem \ref{thm:UBM_inc}, we get the upper bound. Similarly, by picking an \valpha such that $\kappa_{\max}\geq\alpha\lij$ for all $i,j$ we get the lower bound.
\qed\end{proof}

If the saturation rate vector $\vk$ is such that $\kappa_{\max}/\kappa_{\min}$ is constant, without scaling with $n$, it is easy to see that the sum queue lengths under MaxWeight are within a constant factor of the universal lower bound after noting that $\zeta\lij\leq2$ for all $i,j$.

For some values of $\vk$, we get an exact expression for the sum of queue lengths in heavy traffic.
\begin{corollary}\label{cor:nokappa_incomplete}
Suppose that $\vk=\vnu$ in the set of switch systems described in Theorem \ref{thm:UBM_inc}, so that $\vlam^{\epsilon}=\vnu(1-\epsilon)$ for some $\vnu \in \text{Relint}(\cF\linc)$ such that $\numin\triangleq \min\lij \nu\lij >0$ and $0\leq \no,\nt<n$. Then in heavy traffic we have,
\ba
\lim_{\epsilon\to 0} \epsilon \bE\left[\sum_{ij} \bqij\se\right] &= \frac{1}{2}\La \vsig, \vzeta  \Ra\text{ and }\\
\lim_{\epsilon\to 0} \epsilon \bE\left[\bq_{ij}\se\right] &=0 \ \ \forall \ \ i>\no,j>\nt.
\ea
\end{corollary}
\begin{proof}
Since $\ki=1$ and $\ktj=1$ for $i\leq \no$ and $j\leq \nt$, we pick $\alpha\lij=1$ for all $i,j$.
Since the weight vector \valpha satisfies the condition in Theorem 1, the corollary follows. 
\qed
\end{proof}

The following corollary considers the case when exactly one port of the switch is saturated. Under this condition, the switch is said to satisfy the complete resource pooling condition, and was studied in \cite{stolyar2004maxweight}. In this case, the state collapses onto a line. Without loss of generality, the corollary is stated when an input port is saturated.

\begin{corollary}
Consider the set of switch systems operating under the MaxWeight algorithm  as described in Theorem \ref{thm:UBM_inc} with $\no=1$, $\nt=0$ and $\vk=\vnu$ so that $\vlam^{\epsilon}=\vnu(1-\epsilon)$ for some $\vnu \in \text{Relint}(\cF_{10})$ such that $\numin\triangleq \min\lij \nu\lij >0$. Then in heavy traffic we have,
\ba
\lim_{\epsilon\to 0} \epsilon \bE\left[\sum_{j} \bq_{1j}\se\right] &= \frac{\sum_{1j}\sigma^2_{1j}}{2} \text{ and }\\
\lim_{\epsilon\to 0} \epsilon \bE\left[\bq_{ij}\se\right] &=0 \ \ \forall \ \ i>1,j.
\ea
Moreover, MaxWeight is heavy-traffic optimal i.e., $\lim_{\epsilon\to 0} \epsilon \bE\left[\sum_{ij} \bq_{1j}\se\right]$ is minimized under MaxWeight algorithm.
\end{corollary}
The proof follows directly from Theorem \ref{thm:UBM_inc} and the universal lower bound in Proposition \ref{prop:Universal_LB}. In order to clearly see the scaling of queue lengths in terms of $n$, we state the following corollary under Bernoulli arrivals, which again follows from Corollary \ref{cor:nokappa_incomplete} and Proposition \ref{prop:Universal_LB}.

\begin{corollary}\label{cor:bernoulli_incomplete}
Consider the set of switch systems operating under the MaxWeight algorithm, as described in Theorem \ref{thm:UBM_inc}. Suppose that the arrival process for each queue is Bernoulli with the arrival rate vector $\vlam\se$ where $\lambda\se\lij=(1-\epsilon)/2n$ for $i>\no, j>\nt$ and $\lambda\se\lij=(1-\epsilon)/n$ otherwise with $0\leq \no,\nt<n$ so that \no inputs and \nt outputs are saturating. Then, in the heavy traffic limit, we have
\ba
\lim_{\epsilon\to 0} \epsilon \bE\left[\sum_{ij} \bqij\se\right] &= \frac{\no+\nt}{2}\left(1-\frac{1}{n}\right).
\ea
Moreover, under any stable algorithm the sum queue lengths is lower bounded by
\ba
\lim_{\epsilon\to 0} \epsilon \bE\left[\sum_{ij} \bqij\se\right] &\geq  \frac{\max\{\no,\nt\}}{2}\left(1-\frac{1}{n}\right).
\ea
\end{corollary}
So, we know that MaxWeight is within less than a factor of two away from heavy traffic optimality under incomplete saturation. A similar observation was made under the completely saturated case in \cite{switch_arxiv}. It is not clear if the gap is because the lower bound is loose or because MaxWeight is a constant factor away from heavy traffic optimality.
Under the MaxWeight algorithm, the sum of all queue lengths is as if we have $(\no+\nt)$ separate queues, each with variance $(1-\frac{1}{n})$ which is the total variance of the arrivals in each row or column. This is because there are $(\no+\nt)$ independent constraints of the capacity region that are tight in the limit. This is the same reason why the state collapses to the $(\no+\nt)$ dimensional space $\cS\linc$. So, in general the number of tight constraints in the limit is important.

Theorem \ref{thm:UBM_inc} is valid only for incompletely saturated switch $\no<n$ and $\nt<n$. However, a similar result can be proved for the completely saturated case $\no=\nt=n$ as follows.

\begin{corollary}\label{cor:complete_sat}
Consider the set of switch systems operating under the MaxWeight algorithm,  parameterized by  $0<\epsilon<1$ as described in Theorem \ref{thm:UBM_inc}, with the only difference being that, $\no=\nt=n$. Then as long as  $0< \epsilon \leq \numinp/2\|\vk\|$, the steady state queue lengths vector satisfies
\begin{alignat*}{3}
 \left(1-\frac{1}{2n}\right)\frac{\left\|\vsig\se\right\|^2}{\epsilon} - \Biii(\epsilon,n)
 &\leq\bE\left[\La\bvq\se,\valpha \Ra\right]
 &&\leq   \left(1-\frac{1}{2n}\right)\left\|\vsig\se\right\|^2 + \Biv(\epsilon,n)
\end{alignat*}
for any fixed weight vector $\valpha \in \bR^{n^2}$ such that $\La\valpha,\vei \Ra=n\ki$ and $\La\valpha,\vetj \Ra=n\ktj$ for all $1\leq i,j \leq n$, where  $\Biii(\epsilon)$ and $\Biv(\epsilon)$ are $o(\frac{1}{\epsilon})$.
Thus, in the heavy traffic limit as $\epsilon\downarrow 0$, we have
\ba
\lim_{\epsilon\to 0} \epsilon \bE\left[\La\bvq\se,\valpha \Ra\right] &=  \left(1-\frac{1}{2n}\right)\frac{\left\|\vsig\right\|^2}{\epsilon}.
\ea
\end{corollary}
\begin{proof}
Note that Proposition \ref{prop:SSC_inc} is valid in the case when $\no=\nt=n$. Most of the proof of Theorem \ref{thm:UBM_inc} also holds true except for Lemma \ref{lem:projections}. The norm of the projections of unit vectors $\vchi^{(ij)}$ onto the cone $\cS_{nn}$ is different from $\zeta\lij$. This fundamental difference in behavior for the case $\no=\nt=n$ is because of the following reason. For $\no,\nt<n$, the cone $\cS\linc$ is spanned by the vectors $\vei,\vetj$ for $i\leq\no, j\leq\nt$, which are linearly independent and so the cone $\cS\linc$ has dimension $(\no+\nt)$. When $\no=\nt=n$, the vectors $\vei,\vetj$ for $i\leq n, j\leq n$ are not linearly independent because clearly, $\sum_i \vei = \sum_j \vetj$, and so the dimension of cone $\cS_{nn}$ is smaller than $2n$. It can be shown that the cone $\cS_{nn}$ has dimension $(2n-1)$ \cite[page 20]{ziegler1995lectures}.  The proof will be complete once we calculate $\norm{\vchi^{(ij)}\paras}^2$ for all $i,j$. By symmetry, we have that these norms have the same value for all  $i,j$, i.e.,  $\norm{\vchi^{(ij)}\paras}^2 = \xi$ for some $\xi$, which can be calculated as follows. Suppose $\vf_1, \vf_2, ..., \vf_{2n-1}$ is an orthonormal basis of $\cS_{nn}$, we have
\ba
n^2\xi & = \sum_{ij}\norm{\vchi^{(ij)}\paras}^2\\
& = \sum_{ij}\sum_{l=1}^{2n-1}\La\vchi^{(ij)},\vf_l\Ra^2\\
& = \sum_{l=1}^{2n-1}\sum_{ij}\La\vchi^{(ij)},\vf_l\Ra^2\\
& \stackrel{(a)}{=} \sum_{l}^{2n-1}\norm{\vf_l}\paras^2\\
& = 2n-1
\ea
where (a) follows from the fact that $\{\vchi^{(ij)}\}_{ij}$ is an orthonormal basis of $\bR^{(n^2)}$. Replacing $\zeta\lij$ by $\xi=(2n-1)/n^2$, we get the corollary.
\qed\end{proof}
Similar to Corollary \ref{cor:kappa_min}, we can get lower and upper bounds on the sum of queue lengths in heavy traffic using $\kappa_{\max}$ and $\kappa_{\min}$.
We now state the following corollary to illustrate the use of the weight vectors $\valpha$.
\begin{corollary}
Consider the completely saturated switch system in Corollary \ref{cor:complete_sat} with $\vk=\vnu$. Then, in the heavy traffic limit the queue lengths satisfy the following relations.
\ba
\lim_{\epsilon\to 0} \epsilon \bE\left[\sum_{ij} \bqij\se\right] &=  \left(1-\frac{1}{2n}\right)\left\|\vsig\right\|^2\\
\lim_{\epsilon\to 0} \epsilon \bE\left[\sum_{ij} \bqij\nu\lij\se\right] &=  \left(\frac{2n-1}{2n^2}\right) \left\|\vsig\right\|^2\\
\lim_{\epsilon\to 0} \epsilon \bE\left[\sum_{ij} \overline{q}_{i\pi(i)}\se\right] &=  \left(\frac{2n-1}{2n^2}\right) \left\|\vsig\right\|^2 \text{ for any permutation } \pi.
\ea
\end{corollary}
\begin{proof}The proof follows directly from Corollary \ref{cor:complete_sat} by choosing $\valpha=\vone$, $\valpha=n\vnu$ and $\valpha=nP_{\pi}$ respectively, where $P_{\pi}$ is the permutation matrix corresponding to the permutation $\pi$.
\qed
\end{proof}
Note that the first result above is the main result of \cite{switch_arxiv}. Even though the proof of Theorem 1 in \cite{switch_arxiv} is written in different style, it is equivalent to the proof presented here, because the function used there is equivalent to the norm $\norm{\vq\paras}^2$. It can be shown (using orthogonality principle) that for any vector $\vq\in\bR^{n^2}$, its projection onto the space $\cS_{nn}$ is given by
\ba
q_{\parallel\cS ij}=\frac{\sum_i\qij}{n}+\frac{\sum_j\qij}{n}-\frac{\sum\lij\qij}{n^2}.
\ea
Taking the norm, we get
\ba
\norm{\vq\paras}^2 = \frac{1}{n}\left(\sum_i \left(\sum_j\qij \right)^2 +\sum_j \left(\sum_i \qij \right)^2 - \frac{1}{n} \left(\sum\lij \qij \right)^2 \right),
\ea
which is the function used in \cite{switch_arxiv} scaled by $n$. We now present a further special case, which was also studied in \cite{switch_arxiv} to contrast it with the result in Corollary \ref{cor:bernoulli_incomplete}.

\begin{corollary}
Suppose that the traffic of the completely saturated systems in Corollary \ref{cor:complete_sat} is uniform bernoulli with uniform saturation rate, i.e.,  $\vlam\se=(1-\epsilon)/n\vone$, in the heavy traffic we have
\ba
\lim_{\epsilon\to 0} \epsilon \bE\left[\sum_{ij} \bqij\se\right] &= \frac{2n-1}{2}\left(1-\frac{1}{n}\right).
\ea
\end{corollary}
Notice that the behavior of queue length here is similar to $(2n-1)$ separate queues. This is consistent with the discussion after Corollary \ref{cor:bernoulli_incomplete} because in the heavy traffic limit, among the $2n$ constraints in the capacity region, we have only $(2n-1)$ linearly independent ones.


\section{Conclusion \label{sec:conc}}
We consider the heavy-traffic queue length behavior in an  input queued switch operating under the MaxWeight algorithm.
It was recently shown in \cite{switch_arxiv} that, in the heavy-traffic regime, the queue length scales optimally with the size of the switch when all the ports in the switch saturate at the same rate. In this paper, we considered the case when an arbitrary set of ports saturate, and each port is allowed to saturate at different rate. 
%
We obtained an exact heavy traffic characterization of a linear combination of queue lengths and showed that MaxWeight algorithm achieves optimal scaling of sum queue lengths in heavy traffic.
%
%


\begin{acknowledgements}
We thank Prof. Alexander Stolyar and Prof. Michael Harrison for suggesting us the question
about heavy traffic queue length behavior in a switch when different ports saturate at different rates. This research was supported by NSF Grants CNS-1161404, ECCS-1202065 and  ARO Grant W911NF-16-1-0259.
\end{acknowledgements}

\bibliographystyle{spmpsci}      
\bibliography{references}   


\appendix
\section{Proof of Lemma \ref{lem:projections}} \label{appendix:projections}
\begin{proof}
In order to calculate the norm of the projections of the unit vectors $\vchi^{(ij)}$, we will first consider the projections of an arbitrary vector \vx on to the space $\cS\linc$. From the definition of the cone $\cS\linc$, we know that the projection $\vx\paras$ of \vx can be decomposed as
\ba
\vx\paras =\su{i}{1} \wi\vei+\su{j}{2}\wtj\vetj
\ea
From the Orthogonality Principle, we have that
\ba
\La\vx-\vx\paras, \vei\Ra &=0 \quad \text{for } i\leq \no, \\
\La\vx-\vx\paras, \vetj\Ra &=0 \quad \text{for } j\leq \nt,
\ea
Note that
\ba
\La \vei, {\bf e}^{(k)}\Ra &=\begin{cases}
	0 & i\neq k\\
	n & i=k
\end{cases}, \qquad
\La \vetj, \widetilde{{\bf e}}^{(l)}\Ra =\begin{cases}
	0 & j\neq l\\
	n & j=l
\end{cases}\\
\La \vei, \vetj\Ra &=1 \quad \forall i,j.
\ea
Therefore, we get
\ba
n\wi+\su{j}{2}\wtj=\sum_jx\lij \text{ for } i\leq \no  \quad \text{and} \quad
n\wtj+\su{i}{1}\wi=\sum_ix\lij \text{ for } i\leq \nt.
\ea
Defining $\w=\su{i}{1}\wi$ and $\wt=\su{j}{2}\wtj$ and summing each set of equations above, we get
\ba
n\w+n_1\wt=\su{i}{1}\sum_jx\lij  \quad {and} \quad
n\wt+n_2\w=\su{j}{2}\sum_ix\lij
\ea
Solving for $\w$ and $\wt$ we get
\ban
\w = \frac{n\su{i}{1}\sum_jx\lij-n_1\su{j}{2}\sum_ix\lij}{n^2-n_1n_2} \quad \text{and} \quad
\wt = \frac{n\su{j}{2}\sum_ix\lij-n_2\su{i}{1}\sum_jx\lij}{n^2-n_1n_2} \label{eq:W}
\ean
Therefore $\wi$ and $\wtj$ can be written as,
\ban
\wi = \frac{\sum_jx\lij-\wt}{n} \text{ for } i\leq \no   \quad \text{and} \quad
\wtj = \frac{\sum_ix\lij-\w}{n} \text{ for } i\leq \nt   \label{eq:wi}
\ean
The norm of projection $\vx\paras$ is then given by,
\ba
\norm{\vx\paras}^2 &= \La \vx\paras, \vx\paras\Ra\\
&= \La \su{i}{1}\wi\vei+\su{j}{2}\wtj\vetj, \su{i}{1}\wi\vei+\su{j}{2}\wtj\vetj  \Ra\\
&= n\su{i}{1}\wi^2 +\su{i}{1}\wi\su{j}{2}\wtj
    +\su{j}{2}\wtj\su{i}{1}\wi+n\su{j}{2}\wtj^2 \\
&= n\su{i}{1}\wi^2+n\su{j}{2}\wtj^2+2\w\wt \numberthis \label{eq:wnorm}
\ea
We will now consider the four cases in the definition of \vzeta in \eqref{eq:zeta_def}, and calculate $\norm{\vchi^{(ij)}\paras}^2$.
\emph{Case 1: When $1\leq i\leq n_1$ and $1\leq  j\leq n_2$\\}
Using \eqref{eq:W} and \eqref{eq:wi}, it is easy to see that when $\vx=\vchi^{(ij)}$,
	\begin{equation}
	\begin{split}
	\w&=\frac{n-\no}{n^2-\no\nt} \\
	\wk&=\begin{cases}
			\frac{1-\wt}{n} & k=i\\
			\frac{-\wt}{n} & k\neq i
		\end{cases}
	\end{split}
	\qquad
	\begin{split}
	\wt&=\frac{n-\nt}{n^2-\no\nt}\\
	\wtl&=\begin{cases}
			\frac{1-\w}{n} & l=j\\
			\frac{-\w}{n} & l\neq j
		\end{cases}
	\end{split}\nonumber
	\end{equation}

	\ba
	\norm{\vchi^{(ij)}\paras}^2
	&= n\su{k}{1}\wk^2+n\su{l}{2}\wtl^2+2\w\wt\\
	&= n \left[(\no-1)\left(\frac{-\wt}{n}\right)^2+\left(\frac{1-\wt}{n}\right)^2+(\nt-1)\left(\frac{-\w}{n}\right)^2+\left(\frac{1-\w}{n}\right)^2\right] \\
	&\qquad +2\frac{(n-\no)(n-\nt)}{(n^2-\no\nt)^2}\\
	&= \frac{1}{n} \left[\frac{\no(n-\nt)^2}{(n^2-\no\nt)^2}-\frac{2(n-\nt)}{(n^2-\no\nt)}+\frac{\nt(n-\no)^2}{(n^2-\no\nt)^2}-\frac{2(n-\no)}{(n^2-\no\nt)}+2\right] \\
	&\qquad +2\frac{(n-\no)(n-\nt)}{(n^2-\no\nt)^2}\\
&  = \frac{2}{n}-\frac{2n-\no-\nt}{n(n^2-\no\nt)}.
	\ea
	where the last equality is obtained after direct algebraic simplification.
\emph{Case 2: When $1\leq i\leq n_1$ and $j>\nt$\\}
Again, using \eqref{eq:W} and \eqref{eq:wi}, it is easy to see that when $\vx=\vchi^{(ij)}$,
	\begin{equation}
	\begin{split}
	\w&=\frac{n}{n^2-\no\nt} \\
	\wk&=\begin{cases}
			\frac{1-\wt}{n} & k=i\\
			\frac{-\wt}{n} & k\neq i
		\end{cases}
	\end{split}
	\qquad
	\begin{split}
	\wt&=\frac{-\nt}{n^2-\no\nt}\\
	\wtl&=\frac{-\w}{n}
	\end{split}\nonumber
	\end{equation}
	\ba
	\norm{\vchi^{(ij)}\paras}^2
	&= n\su{k}{1}\wk^2+n\su{l}{2}\wtl^2+2\w\wt\\
	&= n \left[(\no-1)\left(\frac{-\wt}{n}\right)^2+\left(\frac{1-\wt}{n}\right)^2+\nt\left(\frac{-\w}{n}\right)^2\right]+2\frac{(n)(-\nt)}{(n^2-\no\nt)^2}\\
	&= \frac{1}{n} \left[\frac{\no(-\nt)^2}{(n^2-\no\nt)^2}-\frac{2(-\nt)}{(n^2-\no\nt)}+\frac{\nt(n)^2}{(n^2-\no\nt)^2}+1\right] -\frac{2n\nt}{(n^2-\no\nt)^2}\\
& = \frac{1}{n}+\frac{\nt}{n(n^2-\no\nt)}
\ea
\emph{Case 3: When $i\geq \no$ and $1\leq j\leq \nt$\\}
Using the same argument as in Case 2, by symmetry, we get
	\ba
	\norm{\vchi^{(ij)}\paras}^2 = \frac{1}{n}+\frac{\no}{n(n^2-\no\nt)}
	\ea
\emph{Case 4: When $i\geq \no$ and $ j\geq \nt$\\}
When $\vx=\vchi^{(ij)}$ with $i\geq \no$ and $ j\geq \nt$ is orthogonal to the space $\cS\linc$ and so, we get $\norm{\vchi^{(ij)}\paras}^2=0$.
\qed
\end{proof}

\end{document}